\documentclass{amsart}
\usepackage{amssymb}
\usepackage{amscd}
\usepackage[latin1]{inputenc}
\usepackage{amsmath}
\usepackage{amsfonts}
\usepackage{latexsym}
\usepackage{verbatim}
\usepackage{graphicx}
\usepackage{epsfig}
\newtheorem{theorem}{Theorem}
\newtheorem{corollary}{Corollary}
\newtheorem{lemma}{Lemma}
\newtheorem{definition}{Definition}
\newtheorem{question}{Question}
\newtheorem{proposition}{Proposition}
\newcommand{\dbar}{d\mkern-6mu\mathchar'26}

\begin{document}
\title[Approximating the hard square entropy constant]{Approximating the hard square entropy constant with probabilistic methods}
\begin{abstract}
For any $\mathbb{Z}^2$ nearest neighbor shift of finite type $X$ and any integer $n \geq 1$, one can define the horizontal strip shift $H_n(X)$ to be the set of configurations on $\mathbb{Z} \times \{1,\ldots,n\}$ which do not contain any forbidden transitions for $X$. Each $H_n(X)$ can be considered as a $\mathbb{Z}$ nearest neighbor shift of finite type, and it is always the case that $\lim_{n \rightarrow \infty} \frac{h^{top}(H_n(X))}{n} = h^{top}(X)$. In this paper, we combine ergodic theoretic techniques with methods from percolation theory and interacting particle systems to show that for the $\mathbb{Z}^2$ hard square shift $\mathcal{H}$, it is in fact the case that $\lim_{n \rightarrow \infty} h^{top}(H_{n+1}(\mathcal{H})) - h^{top}(H_n(\mathcal{H})) = h^{top}(\mathcal{H})$, and that the rate of convergence is at least exponential. A consequence of this is that $h^{top}(\mathcal{H})$ is computable to any tolerance $\frac{1}{n}$ in time polynomial in $n$. We also give an example of a $\mathbb{Z}^2$ block gluing nearest neighbor shift of finite type $Y$ for which $h^{top}(H_{n+1}(Y)) - h^{top}(H_n(Y))$ does not even approach a limit. 
\end{abstract}
\date{}
\author{Ronnie Pavlov}
\address{Ronnie Pavlov\\
Department of Mathematics\\
University of British Columbia\\
1984 Mathematics Road\\
Vancouver, BC V6T 1Z2}
\email{rpavlov@math.ubc.ca}
\urladdr{www.math.ubc.ca/$\sim$rpavlov/}
\thanks{}
\renewcommand{\subjclassname}{MSC 2000}
\maketitle

\section{Introduction}
\label{intro}

Some of the most studied objects in the field of symbolic dynamics are shifts of finite type (or SFTs.) A $\mathbb{Z}^d$ SFT is defined by specifying a finite set $A$, called the alphabet, and a set of forbidden configurations. For any such specification, the associated $\mathbb{Z}^d$ SFT is the set of configurations in $A^{\mathbb{Z}^d}$ in which no forbidden configuration appears. In this paper, we will mostly concern ourselves with nearest neighbor SFTs, which are SFTs for which all forbidden configurations are just pairs of adjacent letters.

To any $\mathbb{Z}^d$ SFT $X$, one can assign a real number $h^{top}(X)$, called its topological entropy. Informally, $h^{top}(X)$ measures the exponential growth rate of the number of configurations which appear in points of $X$. (We postpone a formal definition until Section~\ref{defns}.)
Topological entropy is quite easy to compute for $\mathbb{Z}$ SFTs; to any $\mathbb{Z}$ SFT $X$, one can associate a $0$-$1$ matrix called its transition matrix, and $h^{top}(X)$ is just the logarithm of the Perron eigenvalue of this matrix. For a general introduction to one-dimensional symbolic dynamics and topological entropy, see \cite{LM}.

In general, it is much harder to compute $h^{top}(X)$ for $\mathbb{Z}^2$ SFTs. In fact, there are very few nondegenerate examples of $\mathbb{Z}^2$ SFTs for which the topological entropy has a known closed form. (\cite{B2}, \cite{FiT}, \cite{Kas}, \cite{Lie}) However, one can approximate such a topological entropy by using the easier to compute one-dimensional topological entropies. For any $\mathbb{Z}^2$ nearest neighbor SFT $X$ with alphabet $A$, one can define $H_n(X)$ to be the set of configurations on $\mathbb{Z} \times \{1,\ldots,n\}$ which contain no forbidden pair of adjacent letters. Then $H_n(X)$ can be considered as a $\mathbb{Z}$ nearest neighbor SFT with alphabet the set of legal $n$-high columns in $X$, which we call $A_n(X)$. Two letters $\begin{smallmatrix} a_n\\ \vdots\\ a_1 \end{smallmatrix}$ and $\begin{smallmatrix} b_n\\ \vdots\\ b_1 \end{smallmatrix}$ in $A_n(X)$ may appear consecutively in $H_n(X)$ if and only if $\begin{smallmatrix} a_n b_n\\ \vdots\\ a_1 b_1 \end{smallmatrix}$ is legal in $X$. We can then define $h_n(X) := h^{top}(H_n(X))$, the topological entropy of $H_n(X)$ as a $\mathbb{Z}$ SFT. One can approximate $h^{top}(X)$ via $h_n(X)$; it turns out to be true that $\frac{h_n(X)}{n} \rightarrow h^{top}(X)$ for any $X$. (This is Lemma~\ref{L1} from Section~\ref{prelims}, and we postpone the proof until then.)

One well-studied example of a $\mathbb{Z}^2$ nearest neighbor SFT is the $\mathbb{Z}^2$ hard square shift $\mathcal{H}$, which is the $\mathbb{Z}^2$ nearest neighbor SFT with alphabet $A = \{0,1\}$ where the only forbidden pairs of letters are two adjacent $1$s horizontally or vertically. Since this is the main SFT we study in this paper, we denote $h^{top}(\mathcal{H})$ by $h$, $h_n(\mathcal{H})$ by $h_n$, $H_n(\mathcal{H})$ by $H_n$, and $A_n(\mathcal{H})$ by $A_n$.

There is no known closed form for the topological entropy $h$ of the hard square model, which is also known as the hard square entropy constant. However, there is quite a bit of literature regarding bounds and approximations to $h$. (see \cite{B}, \cite{CW}, \cite{E}, \cite{FJ}) There is, for instance, an algorithm (\cite{Pi}) that lets a computer generate the transition matrix for $H_n$ for any $n$. One can then use these matrices to compute the sequence $h_n$, and use the fact that $\frac{h_n}{n} \rightarrow h$ to approximate $h$. Interestingly, empirical data (\cite{E}, \cite{Pi}) indicates that the differences $h_{n+1} - h_n$ converge much more quickly to $h$; $\frac{h_n}{n}$ seems to converge at a linear rate, whereas $h_{n+1} - h_n$ seems to converge exponentially fast. To our knowledge however, even a proof of the convergence of $h_{n+1} - h_n$ has been an open problem. Our main result shows that this convergence does in fact occur with exponential rate.

\begin{theorem}\label{mainresult}
$\lim_{n \rightarrow \infty} h_{n+1} - h_{n} = h$, and the rate of this convergence is at least exponential.
\end{theorem}

Interestingly, to prove this entirely combinatorial or topological result, we will be using an almost entirely probabilistic or measure-theoretic proof. We use several techniques from the worlds of probability and interacting particle systems, whose definitions and exposition are contained in Section~\ref{prelims}. Our proof relies heavily on some results and techniques from \cite{vdBS}. 

These techniques are quite powerful and have been used to prove results from symbolic dynamics and ergodic theory before; see \cite{BS}, \cite{vdBS}, \cite{Ha}, and \cite{Ha2}. It is our hope that the applications of interacting particle system methods used in this paper will inspire more work on the fascinating interplay between statistical mechanics and symbolic dynamics.

\section{Definitions}
\label{defns}

We here lay out the necessary definitions and terminology for the rest of the paper. An \textbf{alphabet} $A$ will always be a finite set with at least two elements. 

\begin{definition}
The $\mathbb{Z}^d$ \textbf{full shift} on the alphabet $A$ is the set $A^{\mathbb{Z}^d}$. For any full shift $A^{\mathbb{Z}^d}$, we define the $\mathbb{Z}^d$-\textbf{shift action} $\{\sigma_{v}\}_{v \in \mathbb{Z}^d}$ on $A^{\mathbb{Z}^d}$ as follows: for any $v \in \mathbb{Z}^d$ and $x \in A^{\mathbb{Z}^d}$, $(\sigma_{v}(x))(u) = x(v+u)$ for all $u \in \mathbb{Z}^d$.
\end{definition}

\begin{definition}
A $\mathbb{Z}^d$ \textbf{subshift} on an alphabet $A$ is a set $X \subseteq A^{\mathbb{Z}^d}$ with the following two properties:\\

\noindent
{\rm (i)} $X$ is shift-invariant, meaning that for any $x \in X$ and $v \in \mathbb{Z}^d$, 
$\sigma_{v}(x) \in X$.

\noindent 
{\rm (ii)} $X$ is closed in the product topology on $A^{\mathbb{Z}^d}$.
\end{definition}

When the value of $d$ is clear, we will sometimes omit the $\mathbb{Z}^d$ and just use the term subshift.

A \textbf{configuration} $u$ on the alphabet $A$ is any mapping from a non-empty subset $S$ of $\mathbb{Z}^d$ to $A$, where $S$ is called the \textbf{shape} of $u$. 
For any configuration $u$ with shape $S$ and any $T \subseteq S$, denote by $u|_T$ the restriction of $u$ to $T$, i.e. the subconfiguration of $u$ occupying $T$. 

For any integers $a < b$, we use $[a,b]$ to denote $\{a,a+1,\ldots,b\}$.

\begin{definition} 
A $\mathbb{Z}^d$ \textbf{shift of finite type} (or SFT) $X$ is defined by specifying a finite collection $\mathcal{F}$ of finite configurations on $A$, and then defining $X = (A^{\mathbb{Z}^d})_{\mathcal{F}}$ to be the set of $x \in A^{\mathbb{Z}^d}$ such that $x|_S \notin \mathcal{F}$ for all finite $S \subseteq \mathbb{Z}^d$. For any fixed $X$, the \textbf{type} of $X$ is the minimum positive integer $t$ such that for some $\mathcal{F}$ consisting entirely of configurations with shape $[1,t]^d$, $X=(A^{\mathbb{Z}^d})_{\mathcal{F}}$. 
\end{definition}

It is not hard to check that any SFT is a subshift.

Sites $u,v \in \mathbb{Z}^d$ are said to be \textbf{adjacent} if $\sum_{i=1}^d |u_i - v_i| = 1$. If a $\mathbb{Z}^d$ SFT $X$ has forbidden list $\mathcal{F}$ consisting entirely of pairs of adjacent letters, then $X$ is called a $\mathbb{Z}^d$ \textbf{nearest neighbor SFT}. 
In this paper, we will mostly concern ourselves with $d = 1$ or $d = 2$, and all SFTs we consider will be nearest neighbor SFTs.

\begin{definition} 
The $\mathbb{Z}^d$ \textbf{hard square shift} is the nearest neighbor SFT on the alphabet $\{0,1\}$ whose forbidden list $\mathcal{F}$ consists of all pairs of adjacent $1$s in any of the $d$ cardinal directions.
\end{definition}

\begin{definition}
In a nearest neighbor SFT $X$ with alphabet $A$, $a \in A$ is a \textbf{safe symbol} if none of the forbidden configurations in $\mathcal{F}$ contain $A$. In other words, $a$ is a safe symbol if it may legally appear next to any letter of the alphabet in any direction.
\end{definition}

For example, $0$ is a safe symbol for the $\mathbb{Z}^d$ hard square shift.

\begin{definition}
For any $\mathbb{Z}^d$ SFT $X$ with forbidden list $\mathcal{F}$ and any finite configuration $w$ with shape $S$, $w$ is \textbf{locally admissible} in $X$ if $w|_T \notin \mathcal{F}$ for all $T \subseteq S$, and $w$ is \textbf{globally admissible} in $X$ if there exists $x \in X$ for which $x|_S = w$.
\end{definition}

The difference between local and global admissibility is subtle but quite pronounced. It is always quite easy to check whether a configuration is locally admissible, and for $\mathbb{Z}$ SFTs also to check global admissibility. However, for $\mathbb{Z}^2$ SFTs, the question of whether or not a configuration is globally admissible is undecidable. In other words, there does not exist an algorithm which takes as input the set of forbidden configurations $\mathcal{F}$ and a locally admissible configuration $w$, and gives as output an answer to the question of whether $w$ is globally admissible. (\cite{Be}, \cite{Wan}) 

In this paper, we will mostly be concerning ourselves with the $\mathbb{Z}^2$ hard square shift, which we denote by $\mathcal{H}$. All locally admissible configurations in $\mathcal{H}$ are globally admissible, since a locally admissible configuration in $\mathcal{H}$ can always be completed to a point of $\mathcal{H}$ by filling the rest of $\mathbb{Z}^2$ with $0$s. For this reason, we will just refer to any locally admissible or globally admissible configuration in $\mathcal{H}$ as admissible. 

\begin{definition}
The \textbf{language} of a subshift $X$, denoted by $L(X)$, is the set of globally admissible configurations in $X$. The set of configurations with a particular shape $S$ which are in the language of $X$ is denoted by $L_S(X)$. 
\end{definition}

\begin{definition} 
The \textbf{local language} of any $\mathbb{Z}^d$ SFT $X$ with forbidden list $\mathcal{F}$, denoted by $LA(X)$, is the set of all locally admissible configurations in $X$. The set of configurations with shape $S$ which are in the local language of $X$ is denoted by $LA_S(X)$.
\end{definition}

For any configuration $u$ with shape $S$ in $L(X)$, denote by $[u]$ the set $\{x \in X \ : \ x|_S = u\}$, called the \textbf{cylinder set} of $u$. 

\begin{definition}
The \textbf{topological entropy} of a $\mathbb{Z}^d$ subshift $X$, denoted by $h^{top}(X)$, is defined by

\[
h^{top}(X)=\lim_{j_1, j_2, \ldots, j_d \rightarrow \infty} \frac{\ln \big|L_{\prod_{i=1}^d [1,j_i]}(X)\big|}{j_1 j_2 \cdots j_d}.
\]

\end{definition} 

To see why the limit exists, note that the function $f(j_1, \ldots j_d) := \ln |L_{\prod_{i=1}^d [1,j_i]}(X)|$ is subadditive in each coordinate, i.e. for every $i \in [1,d]$ and $a,b > 0$, 
\begin{multline*}
f(j_1, \ldots, j_{i-1}, a+b, j_{i+1}, \ldots, j_d) \leq f(j_1, \ldots, j_{i-1}, a, j_{i+1}, \ldots, j_d) \\ + f(j_1, \ldots, j_{i-1}, b, j_{i+1}, \ldots, j_d).
\end{multline*}

The classical Fekete's subadditivity lemma implies that for any subadditive function $f(n)$ of one variable, $\lim_{n \rightarrow \infty} \frac{f(n)}{n}$ exists. A multivariate version, which can be found in \cite{Cap}, shows that for any function $f(j_1, \ldots, j_d)$ which is subadditive in each variable, 
\[
\lim_{j_1, \ldots, j_d \rightarrow \infty} \frac{f(j_1, \ldots, j_d)}{j_1 j_2 \ldots j_d}
\]
exists (and is invariant of how each $j_i \rightarrow \infty$), and that the limit is equal to the infimum. For $\mathbb{Z}^d$ SFTs, topological entropy may also be computed by using the local language, i.e. if $L_{\prod_{i=1}^d [1,j_i]}(X)$ is replaced by $LA_{\prod_{i=1}^d [1,j_i]}(X)$ in the definition of topological entropy, the limit is unchanged. (\cite{Fr}, \cite{HoMe}) 

\

We will also need some definitions specific to the arguments used in this paper.

We will frequently consider $\mathbb{Z}^d$ as a graph, where two sites are connected by an edge if they are adjacent. For any set $S \subseteq \mathbb{Z}^d$, we identify $S$ with the maximal subgraph of $\mathbb{Z}^d$ with vertex set $S$, i.e. the graph with vertex set $S$ and edges between all pairs of adjacent vertices in $S$.

For any subset $G$ of $\mathbb{Z}^d$, and any set $S \subseteq G$, the \textbf{boundary} of $S$ within $G$, which is denoted by 
$\partial (S,G)$, is the set of $p \in G \setminus S$ which are adjacent to some $q \in S$. If we refer to simply the boundary of a set $S$, or write $\partial S$, then $G$ is assumed to be all of $\mathbb{Z}^d$.

For any integer $i$, we define $R_i = \mathbb{Z} \times \{i\}$, the \textbf{row at height $i$}.

For any partition $\xi$ of a set $S$, and for any $s \in S$, we use $\xi(s)$ to denote the element of $\xi$ in which $s$ lies. If $\xi$ is a partition of the alphabet $A$ of a $\mathbb{Z}^d$ subshift $X$, then $\phi_{\xi}$ is the factor map from $X$ to $\xi^{\mathbb{Z}^d}$ defined by $(\phi_{\xi} x)(v) = \xi(x(v))$ for all $v \in \mathbb{Z}^d$. 

\section{Some preliminaries}
\label{prelims}

We begin by justifying a claim from the introduction.

\begin{lemma}\label{L1}
For any $X$, $\lim_{n \rightarrow \infty} \frac{h_n(X)}{n} = h^{top}(X)$.
\end{lemma}

\noindent
\textit{Proof.} By the comments following the definition of topological entropy, $h^{top}(X) = \inf_{m,n \rightarrow \infty} \frac{\ln |L_{[1,m] \times [1,n]}(X)}{mn}$, and therefore $|L_{[1,m] \times [1,n]}(X)| \geq e^{h^{top}(X)mn}$ for all $m,n \in \mathbb{N}$. By the definition of $h_n(X)$, $h_n(X) = \lim_{m \rightarrow \infty} \frac{\ln |L_{[1,m] \times [1,n]}(X)|}{m}$. Therefore, $\frac{h_n(X)}{n} \geq h^{top}(X)$ for all $n$. Fix any $\epsilon > 0$. By definition of $h^{top}(X)$, there exists $N$ so that for any $m,n > N$, $|L_{[1,m] \times [1,n]}(X)| \leq e^{(h^{top}(X)+\epsilon)mn}$. This means that $\frac{h_n(X)}{n} \leq h^{top}(X) + \epsilon$ for $n > N$. Since $\epsilon$ was arbitrary, we are done.

\begin{flushright}
$\blacksquare$\\
\end{flushright} 

We use several measure-theoretic or probabilistic tools in the proof of Theorem~\ref{mainresult}, chiefly the concepts of percolation, measure-theoretic entropy, stochastic domination, Gibbs measures, and the $\dbar$ metric. We define these notions and state some fundamental theorems relating them in this section. All measures on subshifts considered in this paper are Borel probability measures for the product topology on $A^{\mathbb{Z}^d}$.

We begin by giving a few notations and facts about independent site percolation which will be necessary for our proof. For a detailed introduction to percolation theory, see \cite{Gr}. 

\begin{definition}
For any $0 < p < 1$ and any infinite connected graph $G = (V(G), E(G))$, the \textbf{independent site percolation measure on $G$}, denoted by $P_{p,G}$, is the measure on $\{0,1\}^{V(G)}$ which independently assigns a $1$ with probability $p$ and $0$ with probability $1-p$ at every site in $V(G)$.
\end{definition}

Often a site with a $1$ is said to be open and a site with a $0$ is said to be closed. We define the event $A$ where there exists an infinite connected cluster of $1$s in $G$, and say that $A$ is the event where percolation occurs. One of the foundational principles of percolation theory is that for any countable locally finite graph, there exists a probability $p_c(G)$, called the \textbf{critical probability for site percolation on $G$}, such that for any $p < p_c(G)$, $P_{p,V(G)}(A) = 0$, and for any $p > p_c(G)$, $P_{p,V(G)}(A) > 0$. We most often take $G$ to be the graph representation of $\mathbb{Z}^2$ as described earlier, which is often called the \textbf{square lattice} in the literature. For this reason, the notation $P_p$ with no graph $G$ will always be understood to represent $P_{p,\mathbb{Z}^2}$, and $p_c$ will represent $p_c(\mathbb{Z}^2)$. It was shown in \cite{H} that $p_c > 0.5$, and there have been successive improving lower bounds on $p_c$ since then. (\cite{MeP}, \cite{T}, \cite{vdBE}, \cite{Z}) 

In this paper, we will be concerned only with the case $p < p_c$, where percolation occurs with probability $0$. If $G$ is the square lattice, then this of course implies that $P_p(0 \leftrightarrow \partial([-n,n]^2)) \rightarrow 0$ as $n \rightarrow \infty$, where for any $S \subseteq \mathbb{Z}^2$, $0 \leftrightarrow S$ represents the event where there is a connected path of $1$s starting at $0$ and ending at a point in $S$. In fact, an even stronger statement can be made. The following is a classical theorem from percolation theory, proved by Menshikov. (\cite{M})

\begin{theorem}\label{expdecay}
On the square lattice, for any $p < p_c$, there exist $A$ and $B$ so that $P_p(0 \leftrightarrow \partial([-n,n]^2)) < Ae^{-Bn}$ for all $n$.
\end{theorem}

We now turn to measure-theoretic entropy and conditional measure-theoretic entropy, beginning with finite partitions. For any finite measurable partitions $\xi$ and $\eta$ of a measure space $(X,\mu)$, we make the definitions
\[
H_{\mu}(\xi) = - \sum_{A \in \xi} \mu(A) \log \mu(A) \textrm{ and } H_{\mu}(\xi \ | \ \eta) = - \sum_{A \in \xi, C \in \eta} \mu(A \cap C) \ln \left(\frac{\mu(A \cap C)}{\mu(C)}\right),
\]
where terms with $\mu(A) = 0$ are omitted from the first sum and terms with $\mu(A \cap C) = 0$ are omitted from the second sum.

The following decomposition result can be found in any standard book on ergodic theory (such as \cite{W}):

\begin{theorem}\label{finitedecomp}
For any measure space $(X,\mu)$, any measurable partition $\eta$ of $X$, and any partition $\xi$ of $X$ which is a refinement of $\eta$, $H_{\mu}(\xi) = H_{\mu}(\eta) + H_{\mu}(\xi \ | \ \eta)$.
\end{theorem}

For any measure $\mu$ on a $\mathbb{Z}^d$ subshift which is stationary, i.e. $\mu(B) = \mu(\sigma_{v} B)$ for all $v \in \mathbb{Z}^d$ and measurable $B$, we may define its entropy.

\begin{definition}
For any finite alphabet $A$ and stationary measure $\mu$ on $A^{\mathbb{Z}^d}$, the \textbf{measure-theoretic entropy} of $\mu$ is 
\[
h(\mu)=\lim_{j_1, j_2, \ldots, j_d \rightarrow \infty} \frac{1}{j_1 j_2 \cdots j_d} H_{\mu}\left(\bigvee_{v \in \prod_{i=1}^d [1,j_i]} \sigma_v \mathcal{P}\right),
\]
where $\mathcal{P}$ is the partition of $X$ into cylinder sets determined by the letter at $x(0)$. (i.e. each element of $\mathcal{P}$ is $[a]$ for some $a \in A$.)
\end{definition}

Again, this limit exists (independently of how each $j_i \rightarrow \infty$) and is equal to its infimum by the coordinatewise subadditivity of the function $g(j_1, \ldots, j_d) := H_{\mu}\left(\bigvee_{v \in \prod_{i=1}^d [1,j_i]} \sigma_v \mathcal{P}\right)$ and the already mentioned multivariate generalization of Fekete's subadditivity lemma found in \cite{Cap}.

Alternately, we can write measure-theoretic entropy more concretely:
\[
h(\mu)=\lim_{j_1, j_2, \ldots, j_d \rightarrow \infty} \frac{-1}{j_1 j_2 \cdots j_d} \sum_{w \in A^{\prod_{i=1}^d [1,j_i]}} \mu([w]) \ln \mu([w]),
\]
where terms with $\mu([w]) = 0$ are omitted. 

We will also deal with measure-theoretic conditional entropy. 

\begin{definition}
For any finite alphabet $A$, any stationary measure $\mu$ on $A^{\mathbb{Z}^d}$, and any measurable partition $\xi$ of $A^{\mathbb{Z}^d}$, the \textbf{measure-theoretic conditional entropy} of $\mu$ with respect to $\xi$ is
\[
h(\mu \ | \ \xi) = \lim_{j_1, j_2, \ldots, j_d \rightarrow \infty} \frac{-1}{j_1 j_2 \cdots j_d} H_{\mu}\left(\bigvee_{v \in \prod_{i=1}^d [1,j_i]} \sigma_v \mathcal{P} \ | \ \bigvee_{v \in \prod_{i=1}^d [1,j_i]} \sigma_v \mathcal{\xi}\right),
\]
where again $\mathcal{P}$ represents the partition of $X$ into cylinder sets determined by the letter at $x(0)$.
\end{definition}
Note that when $\xi$ is the partition $\{\varnothing,X\}$, (i.e. $\xi$ ``contains no information'') $h(\mu \ | \ \xi) = h(\mu)$. 

Again there is a more concrete representation for conditional measure-theoretic entropy. We will only deal with the case where $\xi$ is a coarser partition than $\mathcal{P}$, in which case $\xi$ corresponds to some partition of $A$ in an obvious way and we will say $\xi$ was \textbf{induced by} this partition of $A$. For such $\xi$,
\[
h(\mu \ | \ \xi) = \lim_{j_1, j_2, \ldots, j_d \rightarrow \infty} \frac{-1}{j_1 j_2 \cdots j_d} \sum_{w \in A^{\prod_{i=1}^d [1,j_i]}} \mu([w]) \ln\bigg(\frac{\mu([w])}{\mu\left(\left(\bigvee_{v \in \prod_{i=1}^d [1,j_i]} \sigma^v \xi\right)[w]\right)}\bigg),
\]
where terms with $\mu\left(\left(\bigvee_{v \in \prod_{i=1}^d [1,j_i]} \sigma^v \xi\right)[w]\right) = 0$ are omitted.


We note that for any measure $\mu$ on a full shift $A^{\mathbb{Z}^d}$ and any partition $\xi$ of $A^{\mathbb{Z}^d}$ induced by a partition of $A$, the push-forward $\phi_{\xi}(\mu)$ of $\mu$ under the factor map $\phi_{\xi}$ (i.e. the measure $\phi_{\xi}(\mu)$ defined by $(\phi_{\xi}(\mu))(C) = \mu\left(\phi_{\xi}^{-1} C\right)$ for all $C \subset \xi^{\mathbb{Z}^d}$ for which $\phi_{\xi}^{-1} C$ is measurable) is a measure on $\xi^{\mathbb{Z}^d}$. The following proposition follows immediately from Theorem~\ref{finitedecomp} and the definitions of $h(\mu)$, $h(\mu \ | \ \xi)$, and $\phi_{\xi}$: 

\begin{proposition}\label{P2}
For any finite alphabet $A$, any stationary measure $\mu$ on $A^{\mathbb{Z}^d}$, and any partition $\xi$ of $A$,
\[
h(\mu) = h(\mu \ | \ \xi) + h(\phi_{\xi}(\mu)). 
\]
\end{proposition}



Measure-theoretic entropy and topological entropy are related by the following Variational Principle. (See \cite{Mi} for a proof.)

\begin{theorem}\label{VarPri}
For any $\mathbb{Z}^d$ subshift $X$, $h^{top}(X) = \sup h(\mu)$, where $\mu$ ranges over measures whose support is contained in $X$. This supremum is achieved for some such $\mu$.
\end{theorem}

\begin{definition}
A stationary measure $\mu$ supported on a subshift $X$ is called a \textbf{measure of maximal entropy} if $h(\mu) = h^{top}(X)$.
\end{definition}

Measures of maximal entropy will be useful in the proof of Theorem~\ref{mainresult}, since we can rewrite the topological entropies in the theorem as measure-theoretical entropies with respect to measures of maximal entropy. Measures of maximal entropy on nearest neighbor SFTs also have another extremely useful property.

\begin{definition}\label{MRF}
For any finite alphabet $A$ and countable locally finite graph $G = (V(G), E(G))$, a measure $\mu$ on $A^{V(G)}$ is called a $G$-Markov random field (or $G$-MRF) if, for any finite $S \subset V(G)$, any $\eta \in A^S$, any finite $T \subset (V(G) \setminus S)$ s.t. $\partial(S,G) \subseteq T$, and any $\delta \in A^T$ with $\mu([\delta]) \neq 0$,
\[
\mu(x|_S = \eta \ | \ x|_{\partial(S,G)} = \delta|_{\partial(S,G)}) = \mu(x|_S = \eta \ | \ x|_T = \delta).
\]
\end{definition}
Informally, $\mu$ is an MRF if, for any finite $S \subset V(G)$, the sites in $S$ and the sites in $V(G) \setminus (S \cup \partial (S,G))$ are $\mu$-conditionally independent given the sites on $\partial(S,G)$. We note that our definition of MRF differs slightly from the usual one, where the right-hand side would involve conditioning almost surely on an entire configuration on $V(G) \setminus S$ rather than arbitrarily large finite subconfigurations of it. However, the definitions are equivalent and the finite approach leads to simpler calculations and proofs.

\begin{proposition}\label{P1}
{\rm (\cite{BS2}, p. 281, Proposition 1.20)}
For any $\mathbb{Z}^d$ nearest neighbor SFT $X$, all measures of maximal entropy for $X$ are $\mathbb{Z}^d$-MRFs, and for any such measure $\mu$ and any finite shape $S \subseteq \mathbb{Z}^d$, the conditional distribution of $\mu$ on $S$ given any $\delta \in L_{\partial(S,\mathbb{Z}^d)}(X)$ is uniform over all configurations $x \in L_S(X)$ such that the configuration $y$ defined by $y|_S=x$ and $y|_{\partial (S,\mathbb{Z}^d)}=\delta$ is locally admissible in $X$. 
\end{proposition}

In fact we will only use Proposition~\ref{P1} for $d = 1$, where it is a much more classical fact (\cite{Par}), but we state it in full generality here because the conclusion of Proposition~\ref{P1} is related to the well-studied Gibbs measures from statistical physics. In \cite{vdBS}, they study a more general class of measures; in their language, a measure on a $\mathbb{Z}^d$ hard square shift satisfying the conclusion of Proposition~\ref{P1} is called a hard-core measure with all activities $a_i$ equal to $1$. 

\begin{definition}
For any connected subgraph $G$ of the square lattice, $\mu$ is a \textbf{uniform hard-core Gibbs measure} on $G$ if it is a $G$-MRF such that for any finite connected set $B \subset G$ and any admissible $\delta \in L_{\partial(B,G)}(\mathcal{H})$, $\mu(x|_B = \alpha \ | \ x|_{\partial(B,G)} = \delta)$ is uniform over all $\alpha \in A^B$ which are admissible given $\delta$, i.e. the configuration $y$ defined by $y|_B = \alpha$ and $y|_{\partial(B,G)} = \delta$ is in $L(\mathcal{H})$.
\end{definition}

\begin{theorem}\label{Gibbs}
For every infinite connected subgraph $G$ of the square lattice, there is a unique uniform hard-core Gibbs measure on $S$. 
\end{theorem}

\noindent
\textit{Proof.} Theorem 2.3 in \cite{vdBS} implies (in the case where all $a_i = 1$ in their notation) that for any such $G$, there is a unique uniform hard-core Gibbs measure on $G$ if percolation occurs with probability $0$ with respect to $P_{0.5,G}$. We recall that $p_c > 0.5$, and since $G \subseteq \mathbb{Z}^2$, clearly $p_c(G) \geq p_c > 0.5$, and by definition of $p_c(G)$ we are done.

\begin{flushright}
$\blacksquare$\\
\end{flushright}

In fact, we will eventually be able to represent uniform hard-core Gibbs measures on infinite subgraphs $G$ of $\mathbb{Z}^2$ as weak limits of uniform hard-core Gibbs measures on finite $S$, but for this we will need the notion of stochastic dominance. We first need to define the notion of a coupling of a finite set of measures. 

\begin{definition}
For any $n$ and any probability spaces $(X_i, \mu_i)$, $i \in [1,n]$, a \textbf{coupling} of $\mu_1, \mu_2, \ldots, \mu_n$ is a measure $\lambda$ on $\prod_{i=1}^n X_i$ such that for any $j \in [1,n]$ and any $\mu_j$-measurable $B \subseteq X_j$,
\[
\lambda\left(\prod_{i=1}^{j-1} X_i \times B \times \prod_{k=j+1}^n X_k\right) = \mu_j(B).
\] 
\end{definition}

We present two equivalent definitions of stochastic dominance, both of which depend on a partial order $\leq$ on the compact space $A^S$ for some set $S$. We will always assume $\leq$ to be closed, i.e. $\{(x,y) \ : \ x \leq y\} \subset (A^S)^2$ is closed. The equivalence of these definitions is originally due to a result of Strassen (Theorem 11 in \cite{St}, where in his language, $S = T$ and $\epsilon = 0$); also see \cite{Li} for a shorter proof of this equivalence (Theorem 2.4 in \cite{Li}) and a general introduction to interacting particle systems.

\begin{definition}
For any set $S$, any partial ordering $\leq$ on $A^S$, and any measures $\mu$ and $\nu$ on $A^S$, $\mu \leq \nu$ ($\mu$ is \textit{stochastically dominated by} $\nu$ with respect to $\leq$) if there exists a coupling $\lambda$ of $\mu$ and $\nu$ for which $\lambda(\{(x,y) \in (A^S)^2 \ : \ x \leq y\}) = 1$. 
\end{definition}

\begin{definition}
For any set $S$, any partial ordering $\leq$ on $A^S$, and any measures $\mu$ and $\nu$ on $A^S$, $\mu \leq \nu$ ($\mu$ is \textit{stochastically dominated by} $\nu$ with respect to $\leq$) if for any increasing bounded continuous function $f$ from $A^S$ to $\mathbb{R}$ ($f$ is increasing if $f(x) \leq f(y)$ if $x \leq y$), $E_{\mu}(f) \leq E_{\nu}(g)$.
\end{definition}

We will repeatedly make use of three important properties of stochastic dominance.

\begin{lemma}\label{stoch0}
For a partial ordering $\leq$ on $A^S$, define a relation $\leq_T$ on $A^T$ by restricting $\leq$ to $T$. (i.e. $x \leq_T y$ if there exist $x', y' \in A^S$ such that $x'|_T = x$, $y'|_T = y$, and $x' \leq y'$) If $\leq_T$ is a partial order, then for any measures $\mu \leq \nu$ on $A^S$, $\mu|_T \leq_T \nu|_T$.
\end{lemma}

\noindent
\textit{Proof.} This is an obvious consequence of the first definition of stochastic dominance; simply marginalize the coupling $\lambda$ from the first definition of stochastic dominance to get a coupling $\lambda_T$ of $\mu|_T$ and $\nu|_T$ with support contained in $\{(x,y) \in (A^T)^2 \ : \ x \leq_T y\}$.

\begin{flushright}
$\blacksquare$\\
\end{flushright}

\begin{lemma}\label{stoch1}
Stochastic dominance is preserved under weak limits; i.e. if $\mu_n \rightarrow \mu$, $\nu_n \rightarrow \nu$ weakly, and $\mu_n \leq \nu_n$ for all $n$, then $\mu \leq \nu$. 
\end{lemma}

\noindent
\textit{Proof.} This is an obvious consequence of the second definition of stochastic dominance.

\begin{flushright}
$\blacksquare$\\
\end{flushright}

\begin{lemma}\label{stoch2}
If a sequence of measures $\{\mu_n\}$ on $A^S$ is stochastically monotone (i.e. either $\mu_n \leq \mu_{n+1}$ for all $n$ or $\mu_{n+1} \leq \mu_n$ for all $n$), then $\{\mu_n\}$ approaches a weak limit $\mu$.
\end{lemma}

\noindent
\textit{Proof.} We assume that $\{\mu_n\}$ is a stochastically increasing sequence, since the proof is nearly identical for the decreasing case. Since $A^{S}$ is compact, there exists a subsequence of $\mu_n$ which approaches a weak limit. Consider any two subsequences of $\{\mu_n\}$ which each approach weak limits, say $\mu_{n_k} \rightarrow \mu$ and $\mu_{m_k} \rightarrow \mu'$. Then, by passing to subsequences again if necessary, we can assume that $n_1 < m_1 < n_2 < m_2 < \ldots$. Since $\mu_n$ is stochastically increasing, $\mu_{n_i} \leq \mu_{m_i}$ for all $i$ and $\mu_{m_i} \leq \mu_{n_{i+1}}$ for all $i$.

By Lemma~\ref{stoch1}, this means that $\mu' \leq \mu$ and $\mu \leq \mu'$, so $\mu = \mu'$. This means that all weakly convergent subsequences of $\{\mu_n\}$ approach the same limit, and so the sequence itself weakly converges.

\begin{flushright}
$\blacksquare$\\
\end{flushright}

We now define a partial order which is particularly relevant to $\mathcal{H}$. We think of $\mathbb{Z}^2$ as being colored like a checkerboard; $(x,y) \in \mathbb{Z}^2$ is colored black if $x+y$ is even and white if $x+y$ is odd. We define a site-dependent ordering of $\{0,1\}$; for any $v \in \mathbb{Z}^2$, $\preceq_v$ is defined as $0 \preceq_v 1$ if $v$ is black, and $1 \preceq_v 0$ if $v$ is white. We use this site-dependent ordering to define a partial ordering on $\{0,1\}^{S}$ for any $S \subseteq \mathbb{Z}^2$: for any $x,x' \in \{0,1\}^{S}$, $x \preceq x'$ if $x(v) \preceq_v x'(v)$ for all $v \in S$. This in turn defines the stochastic dominance partial ordering on measures on $A^S$ with respect to $\preceq$, which we also denote by $\preceq$.


For any rectangle $R$ and $\delta \in L_{\partial R}(\mathcal{H})$, we define a probability measure $\mu^{\delta}$ on $\{0,1\}^{R}$ which assigns equal probability to all configurations $x$ such that the configuration $y \in \{0,1\}^{R \cup \partial R}$ defined by $y|_R = x$ and $y|_{\partial R} = \delta$ is admissible. (Note that by Proposition~\ref{P1}, $\mu^{\delta}$ is just the conditional probability distribution on $R$, given $\delta$, w.r.t. the measure of maximal entropy $\mu$ for $\mathcal{H}$.) We define a special class of examples: for any $u,d,\ell,r \in \{0,+,-\}$ and any rectangle $R$, define $\delta^{u,d,\ell,r}_R \in L_{\partial R}(\mathcal{H})$ as follows: the symbols $u, d, \ell, r$ determine boundary conditions adjacent to the top, bottom, left, and right edges of $R$. A $+$ means that the sites adjacent to that edge of $R$ are maximal with respect to $\preceq$, i.e. $0$ on white squares and $1$ on black squares. A $-$ means that the sites adjacent to that edge of $R$ are minimal with respect to $\preceq$, i.e. $1$ on white squares and $0$ on black squares. A $0$ means that the sites adjacent to that edge of $R$ are all $0$. We then define $\mu^{u,d,\ell,r}_R$ to be $\mu^{\delta^{u,d,\ell,r}_R}$.

The following theorem states that for the partial order $\preceq$, comparability between two admissible boundary configurations implies stochastic dominance comparability between their associated measures. The theorem is a corollary of Lemma 3.1 from \cite{vdBS}, and the proof is similar to that of Holley's theorem (\cite{Ho}) for the Ising model. 

\begin{theorem}\label{Th1}
For any rectangle $R$ and $\delta, \eta \in L_{\partial R}(\mathcal{H})$ such that $\delta \preceq \eta$, $\mu^{\delta} \preceq \mu^{\eta}$.
\end{theorem}


We can use Theorem~\ref{Th1} to derive stochastic dominance relationships between some of the measures $\mu^{u,d,\ell,r}_{R}$ for different-sized rectangles.

\begin{theorem}\label{Th2}
For any integers $k' \leq k < \ell \leq \ell'$ and $m' \leq m < n$, define rectangles $R = [k,\ell] \times [m,n]$ and $R' = [k', \ell'] \times [m',n]$. Then $\mu^{0,+,+,+}_{R} \succeq \mu^{0,+,+,+}_{R'}|_R$ and $\mu^{0,-,-,-}_{R} \preceq \mu^{0,-,-,-}_{R'}|_R$. 
\end{theorem}

\noindent
\textit{Proof.} We prove only the first inequality, as the second is similar. Our proof mirrors the proof of Proposition 2.5 from \cite{BS}. Since $R \subseteq R'$, we may write $\mu^{0,+,+,+}_{R'}|_R = \mu^{\delta^{0,+,+,+}_{R'}}|_R$ as a weighted average of the measures $\mu^{\eta}$, where $\eta$ ranges over all admissible configurations in $\mathcal{H}$ on $\partial R$ whose top edge is labeled by $0$s. By Theorem~\ref{Th1}, each term in this weighted average is stochastically dominated by $\mu^{0,+,+,+}_{R}$ with respect to $\preceq$, and therefore $\mu^{0,+,+,+}_{R} \succeq \mu^{0,+,+,+}_{R'}|_R$.  

\begin{flushright}
$\blacksquare$\\
\end{flushright}

The proofs of the following two theorems are almost identical.

\begin{theorem}\label{Th2b}
For any integers $k' \leq k < \ell \leq \ell'$ and $m < n \leq n'$, define rectangles $S = [k,\ell] \times [m,n]$ and $S' = [k', \ell'] \times [m,n']$. Then $\mu^{+,0,+,+}_{S} \succeq \mu^{+,0,+,+}_{S'}|_S$ and $\mu^{-,0,-,-}_{S} \preceq \mu^{-,0,-,-}_{S'}|_S$. 
\end{theorem}

\begin{theorem}\label{Th2c}
For any integers $k' \leq k < \ell \leq \ell'$ and $m < n$, define rectangles $T = [k,\ell] \times [m,n]$ and $T' = [k', \ell'] \times [m,n]$. Then $\mu^{0,0,+,+}_{T} \succeq \mu^{0,0,+,+}_{T'}|_T$. 
\end{theorem}

We will also make use of the $\dbar$ topology on probability measures on a full shift $A^{\mathbb{Z}}$. There are many different definitions for the $\dbar$ metric (for a thorough introduction to the subject, see \cite{R}), but the one which we will find most useful is the following.

\begin{definition}
For any stationary measures $\mu$ and $\mu'$ on $A^{\mathbb{Z}}$, 

\[
\dbar(\mu,\mu') = \min_{\lambda \in C(\mu,\mu')} \int d_1(x(0), y(0)) \ d\lambda(x,y),
\]

\noindent
where $C(\mu,\mu')$ is the set of stationary couplings of $\mu$ and $\mu'$ and $d_1$ is the $1$-letter Hamming distance given by $d_1(a,a) = 0$ and $d_1(a,b) = 1$ for $a \neq b$. 
\end{definition}

The $\dbar$ metric is useful for our purposes because of the nice behavior of measure-theoretic entropy in the $\dbar$ topology. We first need a definition:

\begin{definition}
A stationary measure $\mu$ on a $\mathbb{Z}$ subshift $X$ is \textbf{ergodic} if for any shift-invariant measurable set $A \subset X$, i.e. a measurable set $A$ for which $\mu(A \triangle \sigma_n A) = 0$ for all $n$, $\mu(A)$ is $0$ or $1$.
\end{definition}

The following is Theorem 7.9 from \cite{R}.

\begin{theorem}\label{dbarholder}
For any finite alphabet $A$ and ergodic stationary measures $\mu$ and $\nu$ on $A^{\mathbb{Z}}$, if $\dbar(\mu,\nu) = \epsilon$,
then $|h(\mu) - h(\nu)| \leq \epsilon \ln |A| - \epsilon \ln \epsilon - (1 - \epsilon) \ln(1 - \epsilon)$.
\end{theorem}

\section{Main body}
\label{main}

We now restrict our attention to the hard square shift $\mathcal{H}$ and will use our preliminaries to prove some results about measures of maximal entropy on the $\mathbb{Z}$ nearest neighbor shifts of finite type $H_n$. By Theorems~\ref{Th2} and \ref{Th2b}, for any fixed $m \leq n$ and any fixed $K$, the sequences $\left(\mu^{0,+,+,+}_{[-k,k] \times [m,n]}\right)|_{[-K,K] \times [m,n]}$ and $\left(\mu^{+,0,+,+}_{[-k,k] \times [m,n]}\right)|_{[-K,K] \times [m,n]}$ are monotonically decreasing in the stochastic dominance ordering $\preceq$ as $k \rightarrow \infty$. By Lemma~\ref{stoch2}, this implies that for any $K$, these sequences approach weak limits, and so $\mu^{0,+,+,+}_{[-k,k] \times [m,n]}$ and $\mu^{+,0,+,+}_{[-k,k] \times [m,n]}$ approach weak limits, denoted by $\mu^{\stackrel{0}{{\scriptscriptstyle +}}}_{m,n}$ and $\mu^{\stackrel{+}{{\scriptscriptstyle 0}}}_{m,n}$ respectively. An almost identical proof (but with monotonically decreasing marginalizations) shows that $\mu^{0,-,-,-}_{[-k,k] \times [m,n]}$ and $\mu^{-,0,-,-}_{[-k,k] \times [m,n]}$ also approach weak limits as $k \rightarrow \infty$, which we denote by $\mu^{\stackrel{0}{{\scriptscriptstyle -}}}_{m,n}$ and $\mu^{\stackrel{-}{{\scriptscriptstyle 0}}}_{m,n}$ respectively. Finally, by using Theorem~\ref{Th2c} instead of Theorems~\ref{Th2} and \ref{Th2b}, we see that $\mu^{0,0,+,+}_{[-k,k] \times [m,n]}$ approaches a weak limit as $k \rightarrow \infty$, which we denote by $\mu^{\stackrel{0}{{\scriptscriptstyle 0}}}_{m,n}$. (Note: Technically, to discuss weak limits, we need all measures to live on the same space; to deal with this, we could extend each measure to $\{0,1\}^{\mathbb{Z} \times [m,n]}$ by simply appending $0$s to every configuration in the support.) 

\begin{lemma}\label{C1}
For any integer $n$, $\mu^{\stackrel{0}{{\scriptscriptstyle -}}}_{1,n} \preceq \mu^{\stackrel{0}{{\scriptscriptstyle 0}}}_{1,n} \preceq \mu^{\stackrel{0}{{\scriptscriptstyle +}}}_{1,n}$ and $\mu^{\stackrel{-}{{\scriptscriptstyle 0}}}_{1,n} \preceq \mu^{\stackrel{0}{{\scriptscriptstyle 0}}}_{1,n} \preceq \mu^{\stackrel{+}{{\scriptscriptstyle 0}}}_{1,n}$.
\end{lemma}

\noindent
\textit{Proof.} We prove the first set of inequalities only, as the second is similar. For any fixed $k$, $\mu^{0,-,-,-}_{[-k,k] \times [1,n]} \preceq \mu^{0,0,+,+}_{[-k,k] \times [1,n]} \preceq \mu^{0,+,+,+}_{[-k,k] \times [1,n]}$ by Theorem~\ref{Th1}. By letting $k \rightarrow \infty$ and using Lemma~\ref{stoch1}, $\mu^{\stackrel{0}{{\scriptscriptstyle -}}}_{1,n} \preceq \mu^{\stackrel{0}{{\scriptscriptstyle 0}}}_{1,n} \preceq \mu^{\stackrel{0}{{\scriptscriptstyle +}}}_{1,n}$.  

\begin{flushright}
$\blacksquare$\\
\end{flushright}

\begin{lemma}\label{C1.3}
For any integer $n$, $\mu^{\stackrel{0}{{\scriptscriptstyle -}}}_{1,n} \preceq \mu^{\stackrel{0}{{\scriptscriptstyle 0}}}_{0,n}|_{\mathbb{Z} \times [1,n]} \preceq \mu^{\stackrel{0}{{\scriptscriptstyle +}}}_{1,n}$ and $\mu^{\stackrel{-}{{\scriptscriptstyle 0}}}_{1,n} \preceq \mu^{\stackrel{0}{{\scriptscriptstyle 0}}}_{1,n+1}|_{\mathbb{Z} \times [1,n]} \preceq \mu^{\stackrel{+}{{\scriptscriptstyle 0}}}_{1,n}$.
\end{lemma}

\noindent
\textit{Proof.} We again prove the first set of inequalities only, as the second is similar. For any fixed $k$, $\mu^{0,-,-,-}_{[-k,k] \times [1,n]} \preceq \mu^{0,-,-,-}_{[-k,k] \times [0,n]}|_{[-k,k] \times [1,n]} \preceq \mu^{0,0,+,+}_{[-k,k] \times [0,n]}|_{[-k,k] \times [1,n]} \preceq \mu^{0,+,+,+}_{[-k,k] \times [0,n]}|_{[-k,k] \times [1,n]} \preceq \mu^{0,+,+,+}_{[-k,k] \times [1,n]}$ by Theorems~\ref{Th1} and \ref{Th2}. By letting $k \rightarrow \infty$ and using Lemma~\ref{stoch1}, $\mu^{\stackrel{0}{{\scriptscriptstyle -}}}_{1,n} \preceq \mu^{\stackrel{0}{{\scriptscriptstyle -}}}_{0,n}|_{\mathbb{Z} \times [1,n]} \preceq \mu^{\stackrel{0}{{\scriptscriptstyle 0}}}_{0,n}|_{\mathbb{Z} \times [1,n]} \preceq \mu^{\stackrel{0}{{\scriptscriptstyle +}}}_{0,n}|_{\mathbb{Z} \times [1,n]} \preceq \mu^{\stackrel{0}{{\scriptscriptstyle +}}}_{1,n}$, and by removing the second and fourth expressions we are done.

\begin{flushright}
$\blacksquare$\\
\end{flushright}

\begin{theorem}\label{C1.5}
For any $n$, $\mu^{\stackrel{0}{{\scriptscriptstyle 0}}}_{1,n}$ is the unique measure of maximal entropy on $H_n$, and is ergodic.
\end{theorem}

\noindent
\textit{Proof.} By Proposition~\ref{P1}, for any measure $\mu$ of maximal entropy on a $\mathbb{Z}$ nearest neighbor SFT $Y$ and for any $a,b$ letters in the alphabet of $Y$, $\mu(x|_{[m+1,n-1]} = \alpha \ | \ x(m) = a, x(n) = b)$ is uniform over all admissible configurations $\alpha$ given $a$ and $b$. We claim that this implies that any measure of maximal entropy $\mu$ on the $\mathbb{Z}$ nearest neighbor SFT $H_n$, when considered as a measure on $\{0,1\}^{\mathbb{Z} \times [1,n]}$, is a uniform hard-core Gibbs measure on $\mathbb{Z} \times [1,n]$. To see this, consider any finite configurations $w, w' \in \{0,1\}^{S \cup \partial(S,\mathbb{Z} \times [1,n])}$ for some finite $S \subseteq \mathbb{Z} \times [1,n]$ such that $w|_{\partial(S,\mathbb{Z} \times [1,n])} = w'|_{\partial(S,\mathbb{Z} \times [1,n])}$. Then, choose any interval $[l,r]$ so that $S \cup \partial(S,\mathbb{Z} \times [1,n]) \subseteq [l,r] \times [1,n]$, and any configurations $L \in \{0,1\}^{\{l-1\} \times [1,n]}$ and $R \in \{0,1\}^{\{r+1\} \times [1,n]}$ so that $\mu([L] \cap [R]) > 0$. Then, by Proposition~\ref{P1}, all configurations in $L_{[l-1,r+1] \times [1,n]}(\mathcal{H})$ which have $L$ on the left edge and $R$ on the right have the same $\mu$-measure, and so $\mu([w] \cap [L] \cap [R])$ depends only the proportion of such configurations which have restriction $w$ on $S \cup \partial(S,\mathbb{Z} \times [1,n])$. However, since $\mathcal{H}$ is a nearest neighbor SFT, this proportion depends only on the letters on $\partial(S,\mathbb{Z} \times [1,n])$, and so $\mu([w] \cap [L] \cap [R]) = \mu([w'] \cap [L] \cap [R])$. By summing over all such $L,R$, we see that $\mu([w]) = \mu([w'])$, and so $\mu$ is a uniform hard-core Gibbs measure on $\mathbb{Z} \times [1,n]$.

By its definition as a weak limit, it is not hard to check that $\mu^{\stackrel{0}{{\scriptscriptstyle 0}}}_{1,n}$ is also a uniform hard-core Gibbs measure on $\mathbb{Z} \times [1,n]$, and by Theorem~\ref{Gibbs}, there is only one such measure. Therefore, $\mu^{\stackrel{0}{{\scriptscriptstyle 0}}}_{1,n}$ is the unique measure of maximal entropy on $H_n$. It is a standard fact (\cite{Par}) that when a $\mathbb{Z}$ SFT has a unique measure of maximal entropy, it is ergodic.

\begin{flushright}
$\blacksquare$\\
\end{flushright}

We note that by the definition of ergodicity, any marginalization $\mu^{\stackrel{0}{{\scriptscriptstyle 0}}}_{1,n}|_{\bigcup_{a \in A} R_a}$ is also ergodic for $A \subseteq [1,n]$; a shift-invariant set with nontrivial measure for $\mu^{\stackrel{0}{{\scriptscriptstyle 0}}}_{1,n}|_{\bigcup_{a \in A} R_a}$ would yield a shift-invariant set with nontrivial measure for $\mu^{\stackrel{0}{{\scriptscriptstyle 0}}}_{1,n}$.

\begin{theorem}\label{Th9}
For any $k,n$, any even $i \in [1,n]$, and any even $j \in [-k,k]$,

\begin{multline}\notag
0 \leq \mu^{0,-,-,-}_{[-k,k] \times [1,n]}(x(j,i) = 0) - \mu^{0,+,+,+}_{[-k,k] \times [1,n]}(x(j,i) = 0)\\
\leq 2P_{0.5}\Big((j,i) \leftrightarrow \partial\big(([-k,k] \times [1,n]), \mathbb{Z}\times (-\infty,n]\big)\Big) \textrm{ and}
\end{multline}
\begin{multline}\notag
0 \leq \mu^{-,0,-,-}_{[-k,k] \times [1,n]}(x(j,i) = 0) - \mu^{+,0,+,+}_{[-k,k] \times [1,n]}(x(j,i) = 0) \\
\leq 2P_{0.5}\Big((j,i) \leftrightarrow \partial\big(([-k,k] \times [1,n]), \mathbb{Z}\times [1,\infty)\big)\Big).
\end{multline}
The order of the terms in the central differences are reversed when the parity of $i$ or $j$ changes.
\end{theorem}

\noindent
\textit{Proof.} We prove only the first set of inequalities, as the second is completely analogous. For ease of notation, we write $\mu = \mu^{0,-,-,-}_{[-k,k] \times [1,n]}$ and $\mu' = \mu^{0,+,+,+}_{[-k,k] \times [1,n]}$. Since $\mu \preceq \mu'$ by Theorem~\ref{Th1}, and since the function $\chi_{\{x(j,i) = 0\}}$ is a decreasing bounded continuous function on $\mathcal{H}$ with respect to $\preceq$, the inequality $0 \leq \mu(x(j,i) = 0) - \mu'(x(j,i) = 0)$ is clear by the second definition of stochastic dominance.

The second inequality $\mu(x(j,i) = 0) - \mu'(x(j,i) = 0) \leq 2P_{0.5,\mathbb{Z}\times (-\infty,n]}\big((j,i) \leftrightarrow \partial\big(([-k,k] \times [1,n]), \mathbb{Z}\times (-\infty,n]\big)\big)$ will be proved in two steps. We first note that Proposition 3.3 from \cite{vdBS} (where in their notation $\Lambda_n = [-k,k] \times [1,n]$ and the underlying graph $G$ is the subgraph $\mathbb{Z} \times (-\infty,n]$ of the square lattice) implies that
\begin{multline}\label{dbareq2}
\mu(x(j,i) = 0) - \mu'(x(j,i) = 0) \\= (\mu \times \mu')\big(\exists\textrm{path of disagreement from } (j,i) \textrm{ to } \partial\big(([-k,k] \times [1,n]), \mathbb{Z}\times (-\infty,n]\big)\big),
\end{multline}
where a path of disagreement for a pair $(x,y) \in \left( \{0,1\}^{[-k,k] \times [1,n]} \right)^2$ is simply a path of vertices $P$ for which $x(p) \neq y(p)$ for all $p \in P$.

It now suffices to prove that
\begin{multline*}
(\mu \times \mu')\big(\exists\textrm{path of disagreement from } (j,i) \textrm{ to } \\\partial\big(([-k,k] \times [1,n]), \mathbb{Z}\times (-\infty,n]\big)\big) \leq 2P_{0.5}\big((j,i) \leftrightarrow \partial\big(([-k,k] \times [1,n]), \mathbb{Z}\times (-\infty,n]\big)\big).
\end{multline*}

Our proof is just a version of the argument used to prove Corollary 2.2 from \cite{vdBS}, adapted to the finite graph $[-k,k] \times [1,n]$.
We point out first that by the definitions of $\mu$ and $\mu'$, they are MRFs on $[-k,k] \times [1,n]$. The fundamental observation we make is that for any $(j',i') \in [-k,k] \times [1,n]$ and any configurations $\eta, \eta' \in \{0,1\}^{\partial(\{(j',i')\}, [-k,k] \times [1,n])}$, 
\begin{multline*}
(\mu \times \mu')\big(x(j',i') \neq y(j',i') \ : \ x|_{\partial(\{(j',i')\},[-k,k] \times [1,n])} = \eta, \\ y|_{\partial(\{(j',i')\},[-k,k] \times [1,n])} = \eta'\big) \leq 0.5.
\end{multline*}
This is easy to check; from the definitions of $\mu$ and $\mu'$, the conditional distributions 
\[
\mu(x(j',i') \ | \ x|_{\partial(\{(j',i')\},[-k,k] \times [1,n])} = \eta) \textrm{ and } \mu'(y(j',i') \ | \ y|_{\partial(\{(j',i')\},[-k,k] \times [1,n])} = \eta')
\]
are always either uniformly distributed between the letters $0$ and $1$, or entirely concentrated on the letter $0$.

Then we note that since $\mu$ and $\mu'$ are MRFs, for any $(j',i') \in ([-k,k] \times [1,n]) \setminus \{(i,j)\}$ and admissible $\delta, \delta' \in \{0,1\}^{([-k,k] \times [1,n]) \setminus \{(j',i')\}}$, 
\begin{multline}\label{dominance}
(\mu \times \mu')\big((x,y) \textrm{ has a path of disagreement}\\ \textrm{from } (j,i) \textrm{ to } (j',i') \ | \ x|_{([-k,k] \times [1,n]) \setminus \{(j',i')\}} = \delta, y|_{([-k,k] \times [1,n]) \setminus \{(j',i')\}} = \delta'\big) \leq\\
(\mu \times \mu')\big(x(j,i) \neq y(j,i) \ | \ x|_{([-k,k] \times [1,n]) \setminus \{(j',i')\}} = \delta, y|_{([-k,k] \times [1,n]) \setminus \{(j',i')\}} = \delta'\big) \leq 0.5.
\end{multline}
Therefore, the probability measure on $\{0,1\}^{([-k,k] \times [1,n]) \setminus \{(j,i)\}}$ which marks paths of disagreements to $(j,i)$ w.r.t. $\mu \times \mu'$ by $1$s is stochastically dominated by the Bernoulli measure $P_{0.5, ([-k,k] \times [1,n]) \setminus \{(j,i)\}}$ with respect to the standard ordering $\leq$ on $\{0,1\}$, i.e. $0 \leq 1$. More rigorously, if we define a factor map $\tau$ from $\left(\{0,1\}^{([-k,k] \times [1,n]) \setminus \{(j,i)\}}\right)^2$ to $\{0,1\}^{([-k,k] \times [1,n]) \setminus \{(j,i)\}}$ by $(\tau(x,y))(v) = 1$ iff $(x,y)$ has a path of disagreement from $v$ to $(j,i)$, then $\tau(\mu \times \mu') \leq P_{0.5, ([-k,k] \times [1,n]) \setminus \{(j,i)\}}$. (This is proved by constructing a coupling of $\tau(\mu \times \mu')$ and $P_{0.5, ([-k,k] \times [1,n]) \setminus \{(j,i)\}}$ where the Bernoulli trials for $P_{0.5}$ always dominate, which is straightforward by (\ref{dominance}).) Then,

\begin{multline*}\label{dbareq1}
(\mu \times \mu')\left(\exists\textrm{path of disagreement from } (j,i) \textrm{ to } \partial\big(([-k,k] \times [1,n]), \mathbb{Z}\times (-\infty,n]\big)\right)\\
\leq (\mu \times \mu')\big(\exists\textrm{path $\Pi$, not containing $(j,i)$, from a neighbor of $(j,i)$ to}\\
\partial\big(([-k,k] \times [1,n]), \mathbb{Z}\times (-\infty,n]\big) \textrm{ such that for each $p \in \Pi$, there is a path of}\\
\textrm{disagreement from $p$ to $(j,i)$}\big)\\
\leq P_{0.5, ([-k,k] \times [1,n]) \setminus \{(j,i)\}}\big(\exists \textrm{path of $1$s from a neighbor of $(j,i)$ to}\\
\partial\big(([-k,k] \times [1,n]), \mathbb{Z}\times (-\infty,n]\big)\big)
= 2P_{0.5}\big((j,i) \leftrightarrow \partial\big(([-k,k] \times [1,n]), \mathbb{Z}\times (-\infty,n]\big)\big).
\end{multline*}

Combining this with (\ref{dbareq2}) completes the proof. 

\begin{flushright}
$\blacksquare$\\
\end{flushright}

\begin{corollary}\label{C2}
For any $n$, any $i \in [1,n]$, and any $j$,
\[
\bigg|\mu^{\stackrel{0}{{\scriptscriptstyle +}}}_{1,n}(x(j,i) = 0) - \mu^{\stackrel{0}{{\scriptscriptstyle -}}}_{1,n}(x(j,i) = 0)\bigg| \leq 2P_{0.5}\Big((j,i) \leftrightarrow \partial\big((\mathbb{Z} \times [1,n]), \mathbb{Z}\times (-\infty,n]\big)\Big) \textrm{ and}
\]
\[ 
\bigg|\mu^{\stackrel{+}{{\scriptscriptstyle 0}}}_{1,n}(x(j,i) = 0) - \mu^{\stackrel{-}{{\scriptscriptstyle 0}}}_{1,n}(x(j,i) = 0)\bigg| \leq 2P_{0.5}\Big((j,i) \leftrightarrow \partial\big((\mathbb{Z} \times [1,n]), \mathbb{Z}\times [1,\infty)\big)\Big).
\]
\end{corollary}

\noindent
\textit{Proof.} We again prove only the first inequality, as the proof of the second is similar. Let $k \rightarrow \infty$ in Theorem~\ref{Th9} and use the definitions of $\mu^{\stackrel{0}{{\scriptscriptstyle +}}}_{1,n}$ and $\mu^{\stackrel{0}{{\scriptscriptstyle -}}}_{1,n}$ as weak limits. Then note that it is obvious that $P_{0.5,S}(v \leftrightarrow T) \leq P_{0.5}(v \leftrightarrow T)$ for any $S \subset \mathbb{Z}^2$, $v \in \mathbb{Z}^2$, and $T \subseteq S$, since enlarging the universal set $S$ to $\mathbb{Z}^2$ only allows for more possible paths of $1$s from $v$ to $T$.

\begin{flushright}
$\blacksquare$\\
\end{flushright}

Our next result regards closeness of the measures $\mu^{\stackrel{0}{{\scriptscriptstyle 0}}}_{1,n}$ and $\mu^{\stackrel{0}{{\scriptscriptstyle 0}}}_{1,n+1}$ in the $\dbar$ metric when restricted to horizontal strips which are two rows high. We will consider such restrictions as measures on the full shift $(\{0,1\}^{\{0\} \times \{0,1\}})^{\mathbb{Z}}$ for the purposes of the $\dbar$ metric.

\begin{corollary}\label{C3}
For any $n$ and any integer $i \in [1,n-1]$, 
\begin{multline}\notag
\dbar\Big(\mu^{\stackrel{0}{{\scriptscriptstyle 0}}}_{1,n}|_{R_i \cup R_{i+1}}, \mu^{\stackrel{0}{{\scriptscriptstyle 0}}}_{1,n+1}|_{R_i \cup R_{i+1}}\Big) \leq 4P_{0.5}\Big((0,i) \leftrightarrow \partial\big((\mathbb{Z} \times [1,n]), \mathbb{Z}\times [1,\infty)\big)\Big) \\
+ 4P_{0.5}\Big((0,i+1) \leftrightarrow \partial\big((\mathbb{Z} \times [1,n]), \mathbb{Z}\times [1,\infty)\big)\Big) \textrm{ and}\\
\dbar\Big(\mu^{\stackrel{0}{{\scriptscriptstyle 0}}}_{1,n}|_{R_i \cup R_{i+1}}, \mu^{\stackrel{0}{{\scriptscriptstyle 0}}}_{1,n+1}|_{R_{i+1} \cup R_{i+2}}\Big) \leq 4P_{0.5}\Big((0,i) \leftrightarrow \partial\big((\mathbb{Z} \times [1,n]), \mathbb{Z}\times (-\infty,n]\big)\Big) \\
+ 4P_{0.5}\Big((0,i+1) \leftrightarrow \partial\big((\mathbb{Z} \times [1,n]), \mathbb{Z}\times (-\infty,n]\big)\Big).
\end{multline}
\end{corollary}

\noindent
\textit{Proof.} We begin with the first inequality. The proof is fairly similar to that of Lemma 3 from \cite{KKO}, but we cannot apply this directly due to the site-dependence of the ordering $\preceq$. By Lemmas~\ref{C1} and ~\ref{C1.3}, $\mu^{\stackrel{-}{{\scriptscriptstyle 0}}}_{1,n} \preceq \mu^{\stackrel{0}{{\scriptscriptstyle 0}}}_{1,n} \preceq \mu^{\stackrel{+}{{\scriptscriptstyle 0}}}_{1,n}$ and $\mu^{\stackrel{-}{{\scriptscriptstyle 0}}}_{1,n} \preceq \mu^{\stackrel{0}{{\scriptscriptstyle 0}}}_{1,n+1}|_{\mathbb{Z} \times [1,n]} \preceq \mu^{\stackrel{+}{{\scriptscriptstyle 0}}}_{1,n}$. By Lemma~\ref{stoch0}, the same inequalities hold when all four measures are restricted to $R_i \cup R_{i+1}$. Then, by using the first definition of stochastic dominance, the following four couplings exist:

\

A coupling $\lambda_1$ of $\mu^{\stackrel{-}{{\scriptscriptstyle 0}}}_{1,n}|_{R_i \cup R_{i+1}}$ and $\mu^{\stackrel{0}{{\scriptscriptstyle 0}}}_{1,n}|_{R_i \cup R_{i+1}}$ supported on $\{(w,x) \ : \ w \preceq x\}$

A coupling $\lambda_2$ of $\mu^{\stackrel{0}{{\scriptscriptstyle 0}}}_{1,n}|_{R_i \cup R_{i+1}}$ and $\mu^{\stackrel{+}{{\scriptscriptstyle 0}}}_{1,n}|_{R_i \cup R_{i+1}}$ supported on $\{(x,z) \ : \ x \preceq z\}$

A coupling $\lambda_3$ of $\mu^{\stackrel{-}{{\scriptscriptstyle 0}}}_{1,n}|_{R_i \cup R_{i+1}}$ and $\mu^{\stackrel{0}{{\scriptscriptstyle 0}}}_{1,n+1}|_{R_i \cup R_{i+1}}$ supported on $\{(w,y) \ : \ w \preceq y\}$

A coupling $\lambda_4$ of $\mu^{\stackrel{0}{{\scriptscriptstyle 0}}}_{1,n+1}|_{R_i \cup R_{i+1}}$ and $\mu^{\stackrel{+}{{\scriptscriptstyle 0}}}_{1,n}|_{R_i \cup R_{i+1}}$ supported on $\{(y,z) \ : \ y \preceq z\}$

\

By taking the relatively independent coupling of $\lambda_1$ and $\lambda_2$ over the common marginal $\mu^{\stackrel{0}{{\scriptscriptstyle 0}}}_{1,n}|_{R_i \cup R_{i+1}}$, one arrives at a coupling $\lambda_5$ of $\mu^{\stackrel{-}{{\scriptscriptstyle 0}}}_{1,n}|_{R_i \cup R_{i+1}}$, $\mu^{\stackrel{0}{{\scriptscriptstyle 0}}}_{1,n}|_{R_i \cup R_{i+1}}$, and $\mu^{\stackrel{+}{{\scriptscriptstyle 0}}}_{1,n}|_{R_i \cup R_{i+1}}$ supported on $\{(w,x,z) \ : \ w \preceq x \preceq z\}$.

By taking the relatively independent coupling of $\lambda_3$ and $\lambda_4$ over the common marginal $\mu^{\stackrel{0}{{\scriptscriptstyle 0}}}_{1,n+1}|_{R_i \cup R_{i+1}}$, one arrives at a coupling $\lambda_6$ of $\mu^{\stackrel{-}{{\scriptscriptstyle 0}}}_{1,n}|_{R_i \cup R_{i+1}}$, $\mu^{\stackrel{0}{{\scriptscriptstyle 0}}}_{1,n+1}|_{R_i \cup R_{i+1}}$, and $\mu^{\stackrel{+}{{\scriptscriptstyle 0}}}_{1,n}|_{R_i \cup R_{i+1}}$ supported on $\{(w,y,z) \ : \ w \preceq y \preceq z\}$.

\

Finally, by taking the relatively independent coupling of $\lambda_5$ and $\lambda_6$ over the common marginal $\mu^{\stackrel{-}{{\scriptscriptstyle 0}}}_{1,n}|_{R_i \cup R_{i+1}} \times \mu^{\stackrel{+}{{\scriptscriptstyle 0}}}_{1,n}|_{R_i \cup R_{i+1}}$, one arrives at a coupling $\lambda$ of $\mu^{\stackrel{-}{{\scriptscriptstyle 0}}}_{1,n}|_{R_i \cup R_{i+1}}$, $\mu^{\stackrel{0}{{\scriptscriptstyle 0}}}_{1,n}|_{R_i \cup R_{i+1}}$, $\mu^{\stackrel{0}{{\scriptscriptstyle 0}}}_{1,n+1}|_{R_i \cup R_{i+1}}$, and $\mu^{\stackrel{+}{{\scriptscriptstyle 0}}}_{1,n}|_{R_i \cup R_{i+1}}$ supported on 
\[
\{(w,x,y,z) \ : \ w \preceq x \preceq z, w \preceq y \preceq z\} \subset (\{0,1\}^{R_i \cup R_{i+1}})^4.
\]
The measures $\mu^{\stackrel{-}{{\scriptscriptstyle 0}}}_{1,n}$ and $\mu^{\stackrel{+}{{\scriptscriptstyle 0}}}_{1,n}$ are not $\sigma_{(1,0)}$-invariant; in fact their definitions as weak limits imply that $\sigma_{(1,0)}\mu^{\stackrel{-}{{\scriptscriptstyle 0}}}_{1,n} = \mu^{\stackrel{+}{{\scriptscriptstyle 0}}}_{1,n}$. They are, however, $\sigma_{(2,0)}$-invariant, and so we will consider them as measures on $(\{0,1\}^{\{0,1\}^2})^{\mathbb{Z}}$ so that we may treat them as stationary measures. By replacing $\lambda$ by any weak limit of a subsequence of the sequence $\frac{1}{n} \sum_{i=0}^{n-1} \sigma_{(2i,0)} \lambda$, we may also assume that $\lambda$ is $\sigma_{(2,0)}$-invariant. We for now assume that $i$ is even, and claim that
\begin{align*}
\int d_1(w(0), z(0)) \ d\lambda(w,z) & \leq \big(\mu^{\stackrel{-}{{\scriptscriptstyle 0}}}_{1,n}(x(0,i) = 0) - \mu^{\stackrel{+}{{\scriptscriptstyle 0}}}_{1,n}(x(0,i) = 0)\big)\\
& + \big(\mu^{\stackrel{+}{{\scriptscriptstyle 0}}}_{1,n}(x(0,i+1) = 0) - \mu^{\stackrel{-}{{\scriptscriptstyle 0}}}_{1,n}(x(0,i+1) = 0)\big)\\
& + \big(\mu^{\stackrel{+}{{\scriptscriptstyle 0}}}_{1,n}(x(1,i) = 0) - \mu^{\stackrel{-}{{\scriptscriptstyle 0}}}_{1,n}(x(1,i) = 0)\big)\\
& + \big(\mu^{\stackrel{-}{{\scriptscriptstyle 0}}}_{1,n}(x(1,i+1) = 0) - \mu^{\stackrel{+}{{\scriptscriptstyle 0}}}_{1,n}(x(1,i+1) = 0)\big),
\end{align*}
where $w$ and $z$ represent sequences on the alphabet $\{0,1\}^{\{0,1\}^2}$, but $x$ represents a configuration on the alphabet $\{0,1\}$. In fact this is fairly straightforward; we may assume in the integral that $w \preceq z$. This means that $w(0) \neq z(0)$ only when at least one of the inequalities 
\begin{align*}
(w(0))(0,0) < (z(0))(0,0), & \ (w(0))(0,1) > (z(0))(0,1),\\
(w(0))(1,0) > (z(0))(1,0), & \ (w(0))(1,1) < (z(0))(1,1)
\end{align*} 
\noindent
holds. However, since $w$ and $z$ are configurations on $R_i \cup R_{i+1}$, it should be clear that 
\begin{align*}
\lambda(\{(w,z) \ : \ (w(0))(0,0) < (z(0))(0,0)\}) & = (\mu^{\stackrel{-}{{\scriptscriptstyle 0}}}_{1,n}(x(0,i) = 0) - \mu^{\stackrel{+}{{\scriptscriptstyle 0}}}_{1,n}(x(0,i) = 0)\big),\\
\lambda(\{(w,z) \ : \ (w(0))(0,1) > (z(0))(0,1)\}) & = (\mu^{\stackrel{+}{{\scriptscriptstyle 0}}}_{1,n}(x(0,i+1) = 0) - \mu^{\stackrel{-}{{\scriptscriptstyle 0}}}_{1,n}(x(0,i+1) = 0)\big),\\
\lambda(\{(w,z) \ : \ (w(0))(1,0) > (z(0))(1,0)\}) & = (\mu^{\stackrel{+}{{\scriptscriptstyle 0}}}_{1,n}(x(1,i) = 0) - \mu^{\stackrel{-}{{\scriptscriptstyle 0}}}_{1,n}(x(1,i) = 0)\big), \textrm{ and}\\
\lambda(\{(w,z) \ : \ (w(0))(1,1) < (z(0))(1,1)\}) & = (\mu^{\stackrel{-}{{\scriptscriptstyle 0}}}_{1,n}(x(1,i+1) = 0) - \mu^{\stackrel{+}{{\scriptscriptstyle 0}}}_{1,n}(x(1,i+1) = 0)\big).
\end{align*}
Since $\mu^{\stackrel{-}{{\scriptscriptstyle 0}}}_{1,n} = \sigma_{(1,0)} \mu^{\stackrel{+}{{\scriptscriptstyle 0}}}_{1,n}$, the right-hand sides of the first and third inequalities each equal $\big(\mu^{\stackrel{-}{{\scriptscriptstyle 0}}}_{1,n}(x(0,i) = 0) - \mu^{\stackrel{+}{{\scriptscriptstyle 0}}}_{1,n}(x(0,i) = 0)\big)$, and the right-hand sides of the second and fourth inequalities each equal $\big(\mu^{\stackrel{+}{{\scriptscriptstyle 0}}}_{1,n}(x(0,i+1) = 0) - \mu^{\stackrel{-}{{\scriptscriptstyle 0}}}_{1,n}(x(0,i+1) = 0)\big)$. Then, since $\lambda$ is supported on $4$-tuples $(w,x,y,z)$ for which $w \preceq x \preceq z$ and $w \preceq y \preceq z$, and since clearly for such $4$-tuples $w(0) = z(0) \Rightarrow w(0) = x(0) = y(0) = z(0)$, 
\begin{multline}\notag
\dbar\Big(\mu^{\stackrel{0}{{\scriptscriptstyle 0}}}_{1,n}|_{R_i \cup R_{i+1}}, \mu^{\stackrel{0}{{\scriptscriptstyle 0}}}_{1,n+1}|_{R_i \cup R_{i+1}}\Big) \leq \\
\int d_1(x(0), y(0)) \ d\lambda(x,y) \leq \int d_1(w(0), z(0)) \ d\lambda(w,z) \leq\\
2\big(\mu^{\stackrel{-}{{\scriptscriptstyle 0}}}_{1,n}(x(0,i) = 0) - \mu^{\stackrel{+}{{\scriptscriptstyle 0}}}_{1,n}(0,i) = 0)\big) + 2\big(\mu^{\stackrel{+}{{\scriptscriptstyle 0}}}_{1,n}(x(0,i+1) = 0) - \mu^{\stackrel{-}{{\scriptscriptstyle 0}}}_{1,n}(x(0,i+1) = 0)\big),
\end{multline}

\noindent
which by Corollary~\ref{C2} is bounded from above by 
\[
4P_{0.5}\Big((0,i) \leftrightarrow \partial\big((\mathbb{Z} \times [1,n]), \mathbb{Z}\times [1,\infty)\big)\Big) + 4P_{0.5}\Big((0,i+1) \leftrightarrow \partial\big((\mathbb{Z} \times [1,n]), \mathbb{Z}\times [1,\infty)\big)\Big).
\]
Due to the earlier rescaling forced by the non-shift invariance of $\mu^{\stackrel{-}{{\scriptscriptstyle 0}}}_{1,n}$ and $\mu^{\stackrel{+}{{\scriptscriptstyle 0}}}_{1,n}$, this is in fact a bound on the $\dbar$ distance between $\mu^{\stackrel{0}{{\scriptscriptstyle 0}}}_{1,n}|_{R_i \cup R_{i+1}}$ and $\mu^{\stackrel{0}{{\scriptscriptstyle 0}}}_{1,n+1}|_{R_i \cup R_{i+1}}$ as measures on $(\{0,1\}^{\{0,1\}^2})^{\mathbb{Z}}$ rather than $(\{0,1\}^{\{0\} \times \{0,1\}})^{\mathbb{Z}}$, but clearly the $\dbar$ distance in the latter case is even smaller.

To prove the second inequality, simply use $\mu^{\stackrel{0}{{\scriptscriptstyle 0}}}_{0,n}$ instead of $\mu^{\stackrel{0}{{\scriptscriptstyle 0}}}_{1,n+1}$, and note that $\mu^{\stackrel{0}{{\scriptscriptstyle 0}}}_{0,n}|_{R_i \cup R_{i+1}} = \mu^{\stackrel{0}{{\scriptscriptstyle 0}}}_{1,n+1}|_{R_{i+1} \cup R_{i+2}}$.

The proofs when $i$ are odd are almost identical, except that the orders of all differences above need to be switched, which does not affect the final inequality.

\begin{flushright}
$\blacksquare$\\
\end{flushright}

Since $0.5 < p_c$ on the square lattice, the following is clear from Theorem~\ref{expdecay} and Corollary~\ref{C3}.

\begin{theorem}\label{closerows}
There exist $A,B > 0$ so that for any $n$ and $i \in [1,n-1]$, 
\[
\dbar\Big(\mu^{\stackrel{0}{{\scriptscriptstyle 0}}}_{1,n}|_{R_i \cup R_{i+1}}, \mu^{\stackrel{0}{{\scriptscriptstyle 0}}}_{1,n+1}|_{R_i \cup R_{i+1}}\Big) \leq Ae^{-B(n-i)} \textrm{ and}
\]
\[
\dbar\Big(\mu^{\stackrel{0}{{\scriptscriptstyle 0}}}_{1,n}|_{R_i \cup R_{i+1}}, \mu^{\stackrel{0}{{\scriptscriptstyle 0}}}_{1,n+1}|_{R_{i+1} \cup R_{i+2}}\Big) \leq Ae^{-Bi}.
\]
\end{theorem}

We note that clearly Theorem~\ref{closerows} also implies that 
\[
\dbar\Big(\mu^{\stackrel{0}{{\scriptscriptstyle 0}}}_{1,n}|_{R_i}, \mu^{\stackrel{0}{{\scriptscriptstyle 0}}}_{1,n+1}|_{R_i}\Big) \leq Ae^{-B(n-i)} \textrm{ and } \dbar\Big(\mu^{\stackrel{0}{{\scriptscriptstyle 0}}}_{1,n}|_{R_i}, \mu^{\stackrel{0}{{\scriptscriptstyle 0}}}_{1,n+1}|_{R_{i+1}}\Big) \leq Ae^{-Bi};
\]
\noindent
either inequality can be proved by considering a restriction of the coupling $\lambda$ that achieves the analogous $\dbar$ distance in Theorem~\ref{closerows} and noting that restricting from a strip two rows high to a single row cannot introduce new disagreements.

Now, we can prove Theorem~\ref{mainresult} by using measure-theoretic conditional entropies. We first need some notation and a preliminary theorem. For any $H_n$, any stationary measure $\mu$ on $H_n$, and any adjacent intervals $I, J \subseteq [1,n]$, we partition the alphabet $A_n = L_{\{0\} \times [1,n]}(\mathcal{H})$ of $H_n$ by the letters on $I \cup J$, and call this partition $\xi_{I \cup J}$. We also partition $A_n$ by the letters on $I$, and call this partition $\xi_I$. Then we make the notations
\[
h_{\mu}\Big(\bigcup_{i \in I} R_i\Big) := h\left(\phi_{\xi_I}(\mu)\right) \textrm{ and}
\]
\[
h_{\mu}\Big(\bigcup_{j \in J} R_j \ | \ \bigcup_{i \in I} R_i\Big) := h\left(\phi_{\xi_{I \cup J}}(\mu) \ | \ \xi_I\right).
\] 

(For the sake of completeness, we note that for any $I$, $\phi_{\xi_I}(\mu)$ is essentially just $\mu|_{\bigcup_{i \in I} R_i}$; we use the partition notation to more easily apply Proposition~\ref{P2}.) We note that $h_{\mu}\Big(\bigcup_{i \in I} R_i\Big)$ can also be thought of as $h(\mu|_{\bigcup_{i \in I} R_i})$. We also note that by Proposition~\ref{P2}, for any $I$ and $J$,
\[
h_{\mu}\Big(\bigcup_{k \in I \cup J} R_k\Big) = h_{\mu}\Big(\bigcup_{i \in I} R_i\Big) + h_{\mu}\Big(\bigcup_{j \in J} R_j \ | \ \bigcup_{i \in I} R_i\Big).
\]
For uniform hard-core Gibbs measures on $\mathbb{Z} \times [1,n]$, we will prove an important fact about these conditional measure-theoretic entropies, which can be thought of as a two-dimensional entropic analogue of the fact that the future and past of a one-dimensional Markov chain are conditionally independent given the present.

\begin{theorem}\label{gibbsindep}
For any $n$ and any adjacent intervals $I, J \subseteq [1,n]$,
\[
h_{\mu^{\stackrel{0}{{\scriptscriptstyle 0}}}_{1,n}}\Big(\bigcup_{j \in J} R_j \ | \ \bigcup_{i \in I} R_i\Big) = h_{\mu^{\stackrel{0}{{\scriptscriptstyle 0}}}_{1,n}}\Big(\bigcup_{j \in J} R_j \ | \ R_i\Big),
\]
\noindent
where $i \in I$ is the element of $I$ adjacent to $J$.
\end{theorem}

\noindent
\textit{Proof.} We will only prove the theorem for the case where $J$ is above $I$, i.e. $I = [i',i]$ and $J=[i+1,j]$, as the other case is trivially similar. Also, for this proof, given a finite set of configurations $\alpha_i \in A^{S_i}$, $1 \leq i \leq k$, for which the shapes $S_i$ are pairwise disjoint, we denote by $\alpha_1 \alpha_2 \ldots \alpha_k$ the concatenation of the $\alpha_i$, i.e. the configuration on $\bigcup_{i=1}^k S_i$ for which $(\alpha_1 \alpha_2 \ldots \alpha_k)|_{S_i} = \alpha_i$ for $1 \leq i \leq k$.

For readability, we abbreviate $\mu^{\stackrel{0}{{\scriptscriptstyle 0}}}_{1,n}$ by $\mu$ in this proof. By definition, since the support of $\mu$ is contained in $H_n$,
\begin{align*}
h_{\mu}\Big(\bigcup_{j \in J} R_j \ | \ \bigcup_{i \in I} R_i\Big) & = \lim_{k \rightarrow \infty} \frac{1}{2k+1} \sum_{\substack{w \in L_{[-k,k] \times I}(\mathcal{H}),\\x \in L_{[-k,k] \times J}(\mathcal{H})}} \mu([w] \cap [x]) \ln \Big(\frac{\mu([w])}{\mu([w] \cap [x])}\Big)\\
& = \lim_{k \rightarrow \infty} \frac{1}{2k+1} \int_{H_n} \ln \Big(\frac{\mu([w])}{\mu([w] \cap [x])}\Big) \ d\mu(w,x).
\end{align*}
We make the decomposition 
\begin{equation}\label{*}
\mu([w] \cap [x]) = \sum_{\substack{L \in L_{\{-k-1\} \times [1,n]}(\mathcal{H}),\\R \in L_{\{k+1\} \times [1,n]}(\mathcal{H})}} \mu([w] \cap [x] \cap [L] \cap [R]).
\end{equation}

We recall from its definition as a weak limit that $\mu$ is a uniform hard-core Gibbs measure on $\mathbb{Z} \times [1,n]$, and so for any such $L$ and $R$,
\[
\mu([w] \cap [x] \cap [L] \cap [R]) = \mu([L] \cap [R]) \frac{|\{u \in L_{[-k,k] \times ([1,i'-1] \cup [j+1,n])}(\mathcal{H}) \ : \ LuwxR \in L(X)\}|}{|\{u \in L_{[-k,k] \times ([1,n])}(\mathcal{H}) \ : \ LuR \in L(\mathcal{H})\}|}.
\]
\begin{figure}[h]
\centering
\includegraphics[scale=0.8]{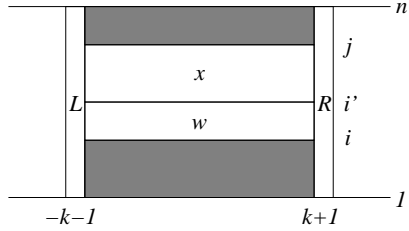}
\caption{$w$, $x$, $L$, and $R$}
\label{F6}
\end{figure}

In Figure~\ref{F6}, $\{u \in L_{[-k,k] \times ([1,i'-1] \cup [j+1,n])}(\mathcal{H}) \ : \ LuwxR \in L(X)\}$ is the set of configurations which can legally fill the shaded area.

We may similarly decompose $\mu([w])$:
\begin{equation}\label{**}
\mu([w]) = \sum_{\substack{L \in L_{\{-k-1\} \times [1,n]}(\mathcal{H}),\\R \in L_{\{k+1\} \times [1,n]}(\mathcal{H})}} \mu([w] \cap [L] \cap [R]).
\end{equation}
Since $\mu$ is a uniform hard-core Gibbs measure on $\mathbb{Z} \times [1,n]$, for any such $L$ and $R$
\[
\mu([w] \cap [L] \cap [R]) = \mu([L] \cap [R]) \frac{|\{u \in L_{[-k,k] \times ([1,i'-1] \cup [i+1,n])}(\mathcal{H}) \ : \ LuwR \in L(X)\}|}{|\{u \in L_{[-k,k] \times ([1,n])}(\mathcal{H}) \ : \ LuR \in L(\mathcal{H})\}|}.
\]
By (\ref{**}) and (\ref{***}), for any $L$ and $R$ such that $\mu([L] \cap [R]) > 0$,
\begin{equation}\label{***}
\frac{\mu([w] \cap [L] \cap [R])}{\mu([w] \cap [x] \cap [L] \cap [R])} = \frac{|\{u \in L_{[-k,k] \times ([1,i'-1] \cup [i+1,n])}(\mathcal{H}) \ : \ LuwR \in L(X)\}|}{|\{u \in L_{[-k,k] \times ([1,i'-1] \cup [j+1,n])}(\mathcal{H}) \ : \ LuwxR \in L(X)\}|}.
\end{equation}
Since $\mathcal{H}$ is a nearest neighbor SFT, 
\begin{multline*}
|\{u \in L_{[-k,k] \times ([1,i'-1] \cup [i+1,n])}(\mathcal{H}) \ : \ LuwR \in L(X)\}| = \\ \left(|\{u' \in L_{[-k,k] \times [1,i'-1]}(\mathcal{H}) \ : \ Lu'wR \in L(X)\}| \right)\\ \cdot \left( |\{u'' \in L_{[-k,k] \times [i+1,n]}(\mathcal{H}) \ : \ Lu''wR \in L(X)\}| \right)
\end{multline*}
and
\begin{multline*}
|\{u \in L_{[-k,k] \times ([1,i'-1] \cup [j+1,n])}(\mathcal{H}) \ : \ LuwxR \in L(X)\}| = \\ \left(|\{u' \in L_{[-k,k] \times [1,i'-1]}(\mathcal{H}) \ : \ Lu'wR \in L(X)\}| \right)\\ \cdot \left( |\{u'' \in L_{[-k,k] \times [j+1,n]}(\mathcal{H}) \ : \ Lu''xR \in L(X)\}| \right).
\end{multline*}

Therefore, (\ref{***}) implies
\begin{equation}\label{****}
\frac{\mu([w] \cap [L] \cap [R])}{\mu([w] \cap [x] \cap [L] \cap [R])} = \frac{|\{u \in L_{[-k,k] \times [i+1,n]}(\mathcal{H}) \ : \ LuwR \in L(X)\}|}{|\{u \in L_{[-k,k] \times [j+1,n]}(\mathcal{H}) \ : \ LuxR \in L(X)\}|}.
\end{equation}
\begin{figure}[h]
\centering
\includegraphics[scale=0.8]{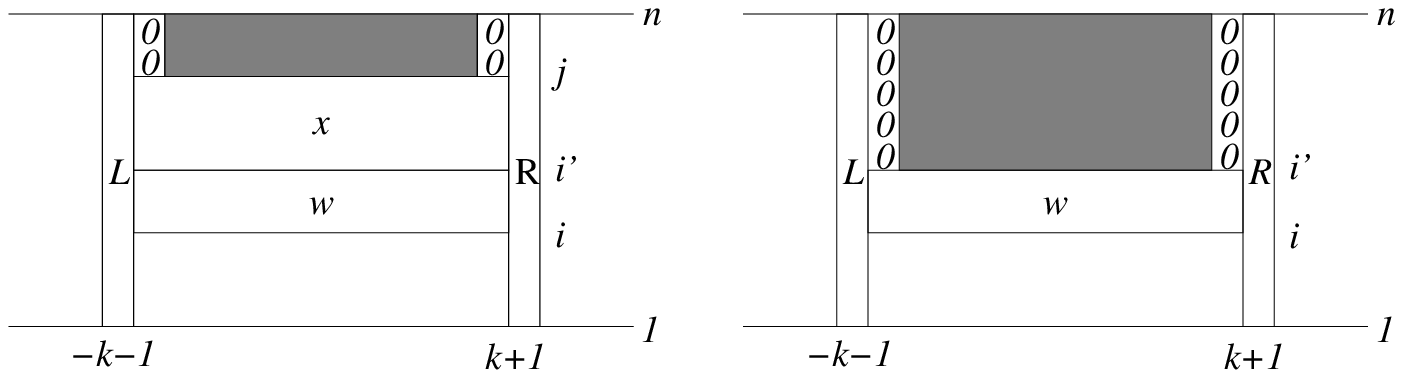}
\caption{}
\label{F7}
\end{figure}

In Figure~\ref{F7}, we see that since $0$ is a safe symbol of $\mathcal{H}$, any configuration $u \in L_{[-k+1,k-1] \times [i+1,n]}(\mathcal{H})$ for which $uw \in L(\mathcal{H})$ may be extended in at least one way to a configuration $u' \in L_{[-k,k] \times [i+1,n]}(\mathcal{H})$ for which $Lu'wR \in L(\mathcal{H})$ (by placing columns of $0$s to the left and right), and any configuration $u \in L_{[-k+1,k-1] \times [j+1,n]}(\mathcal{H})$ for which $ux \in L(\mathcal{H})$ may be extended in at least one way to a configuration $u' \in L_{[-k,k] \times [j+1,n]}(\mathcal{H})$ for which $Lu'xR \in L(\mathcal{H})$ (by placing columns of $0$s to the left and right.) Also, clearly there are at most $2^{2(n-i)}$ possible such extensions in the first case and at most $2^{2(n-j)}$ possible such extensions in the second. Therefore, 
\begin{multline*}
|\{u \in L_{[-k+1,k-1] \times [j+1,n]}(\mathcal{H}) \ : \ ux \in L(X)\}| \\ \leq |\{u \in L_{[-k,k] \times [j+1,n]}(\mathcal{H}) \ : \ LuxR \in L(X)\}|\\ \leq 2^{2(n-j)} |\{u \in L_{[-k+1,k-1] \times [j+1,n]}(\mathcal{H}) \ : \ ux \in L(X)\}| \textrm{ and}
\end{multline*}
\begin{multline*}
|\{u \in L_{[-k+1,k-1] \times [i+1,n]}(\mathcal{H}) \ : \ uw \in L(X)\}| \\ \leq |\{u \in L_{[-k,k] \times [i+1,n]}(\mathcal{H}) \ : \ LuwR \in L(X)\}|\\ \leq 2^{2(n-i)} |\{u \in L_{[-k+1,k-1] \times [i+1,n]}(\mathcal{H}) \ : \ uw \in L(X)\}|.
\end{multline*}
Then by (\ref{****}), for any choices of $L,R$ for which $\mu([L] \cap [R]) > 0$,
\[
\frac{\frac{\mu([w] \cap [L] \cap [R])}{\mu([w] \cap [x] \cap [L] \cap [R])}}{\frac{|\{u \in L_{[-k+1,k-1] \times [i+1,n]}(\mathcal{H}) \ : \ uw \in L(X)\}|}{|\{u \in L_{[-k+1,k-1] \times [j+1,n]}(\mathcal{H}) \ : \ ux \in L(X)\}|}} \in \big[2^{-2(n-i)}, 2^{2(n-j)}\big], \textrm{ and so by (\ref{*}) and (\ref{**}),}
\]
\[
\frac{\frac{\mu([w])}{\mu([w] \cap [x])}}{\frac{|\{u \in L_{[-k+1,k-1] \times [i+1,n]}(\mathcal{H}) \ : \ uw \in L(X)\}|}{|\{u \in L_{[-k+1,k-1] \times [j+1,n]}(\mathcal{H}) \ : \ ux \in L(X)\}|}} \in \big[2^{-2(n-i)}, 2^{2(n-j)}\big].
\]
The original conditional entropy $h_{\mu}(\bigcup_{j \in J} R_j \ | \ \bigcup_{i \in I} R_i)$ is
\[
\lim_{k \rightarrow \infty} \frac{1}{2k+1} \int_{H_n} \ln \Big(\frac{\mu([w])}{\mu([w] \cap [x])}\Big) \ d\mu(w,x), \textrm{ which is equal to}
\]
\[
\lim_{k \rightarrow \infty} \frac{1}{2k+1} \int_{H_n} \ln \left(\frac{|\{u \in L_{[-k+1,k-1] \times [i+1,n]}(\mathcal{H}) \ : \ uw \in L(X)\}|}{|\{u \in L_{[-k+1,k-1] \times [j+1,n]}(\mathcal{H}) \ : \ ux \in L(X)\}|}\right) \ d\mu(w,x)
\]
\noindent
since the difference between the functions inside the integrals is bounded as $k \rightarrow \infty$.

We now note that this expression does not depend on the left endpoint $i'$ of $I$, and so
\[
h_{\mu}\Big(\bigcup_{j \in J} R_j \ | \ \bigcup_{i \in I} R_i\Big) = h_{\mu}\Big(\bigcup_{j \in J} R_j \ | \ R_i\Big).
\]

\begin{flushright}
$\blacksquare$\\
\end{flushright}

\noindent
\textit{Proof of Theorem~\ref{mainresult}.} By Theorem~\ref{C1.5}, $h_{n+1} = h\left(\mu^{\stackrel{0}{{\scriptscriptstyle 0}}}_{1,n+1}\right)$ and $h_n = h\left(\mu^{\stackrel{0}{{\scriptscriptstyle 0}}}_{1,n}\right)$. By using Proposition~\ref{P2}, we may decompose these entropies as follows:
\begin{align*}
h_n = h(\mu^{\stackrel{0}{{\scriptscriptstyle 0}}}_{1,n}) 
& \ = \ h_{\mu^{\stackrel{0}{{\scriptscriptstyle 0}}}_{1,n}}\big(R_{\lfloor \frac{n}{2} \rfloor}\big) \\ 
& \ + \ \sum_{j=\lfloor \frac{n}{2} \rfloor+1}^{n} h_{\mu^{\stackrel{0}{{\scriptscriptstyle 0}}}_{1,n}}\bigg(R_j \ | \ \bigcup_{i=\lfloor \frac{n}{2} \rfloor}^{j-1} R_i\bigg) \\ 
& \ + \ \sum_{k=1}^{\lfloor \frac{n}{2} \rfloor-1} h_{\mu^{\stackrel{0}{{\scriptscriptstyle 0}}}_{1,n}}\bigg(R_k \ | \ \bigcup_{i=k+1}^{n} R_i\bigg) \textrm{ and}\\
h_{n+1} = h(\mu^{\stackrel{0}{{\scriptscriptstyle 0}}}_{1,n+1})
& \ = \ h_{\mu^{\stackrel{0}{{\scriptscriptstyle 0}}}_{1,n+1}}\big(R_{\lfloor \frac{n}{2} \rfloor}\big) \\ 
& \ + \ h_{\mu^{\stackrel{0}{{\scriptscriptstyle 0}}}_{1,n+1}}\big(R_{\lfloor \frac{n}{2} \rfloor+1} \ | \ R_{\lfloor \frac{n}{2} \rfloor}\big) \\
& \ + \ \sum_{j=\lfloor \frac{n}{2} \rfloor+2}^{n+1} h_{\mu^{\stackrel{0}{{\scriptscriptstyle 0}}}_{1,n+1}}\bigg(R_j \ | \ \bigcup_{i=\lfloor \frac{n}{2} \rfloor}^{j-1} R_i\bigg) \\
& \ + \ \sum_{k=1}^{\lfloor \frac{n}{2} \rfloor-1} h_{\mu^{\stackrel{0}{{\scriptscriptstyle 0}}}_{1,n+1}}\bigg(R_k \ | \ \bigcup_{i=k+1}^{n+1} R_i\bigg).
\end{align*}
By Theorem~\ref{gibbsindep}, these decompositions may be rewritten as
\begin{align*}
h_n = h(\mu^{\stackrel{0}{{\scriptscriptstyle 0}}}_{1,n})
& \ = \ h_{\mu^{\stackrel{0}{{\scriptscriptstyle 0}}}_{1,n}}\big(R_{\lfloor \frac{n}{2} \rfloor}\big) \\
& \ + \ \sum_{j=\lfloor \frac{n}{2} \rfloor+1}^{n} h_{\mu^{\stackrel{0}{{\scriptscriptstyle 0}}}_{1,n}}(R_j \ | \ R_{j-1}) \\ \
& \ + \ \sum_{k=1}^{\lfloor \frac{n}{2} \rfloor-1} h_{\mu^{\stackrel{0}{{\scriptscriptstyle 0}}}_{1,n}}(R_k \ | \ R_{k+1}) \textrm{ and}
\end{align*}
\begin{align*}
h_{n+1} = h(\mu^{\stackrel{0}{{\scriptscriptstyle 0}}}_{1,n+1})
& \ = \ h_{\mu^{\stackrel{0}{{\scriptscriptstyle 0}}}_{1,n+1}}\big(R_{\lfloor \frac{n}{2} \rfloor}\big) \\
& \ + \ h_{\mu^{\stackrel{0}{{\scriptscriptstyle 0}}}_{1,n+1}}\big(R_{\lfloor \frac{n}{2} \rfloor+1} \ | \ R_{\lfloor \frac{n}{2} \rfloor}\big) \\ 
& \ + \ \sum_{j=\lfloor \frac{n}{2} \rfloor+2}^{n+1} h_{\mu^{\stackrel{0}{{\scriptscriptstyle 0}}}_{1,n+1}}(R_j \ | \ R_{j-1}) \\ 
& \ + \ \sum_{k=1}^{\lfloor \frac{n}{2} \rfloor-1} h_{\mu^{\stackrel{0}{{\scriptscriptstyle 0}}}_{1,n+1}}(R_k \ | \ R_{k+1}).
\end{align*}
By taking the difference, we see that $h_{n+1} - h_n = $
\begin{align}
\label{line1} h_{\mu^{\stackrel{0}{{\scriptscriptstyle 0}}}_{1,n+1}}\big(R_{\lfloor \frac{n}{2} \rfloor}\big) & \ - \ h_{\mu^{\stackrel{0}{{\scriptscriptstyle 0}}}_{1,n}}\big(R_{\lfloor \frac{n}{2} \rfloor}\big)\\
\label{line2}+ \ \sum_{j=\lfloor \frac{n}{2} \rfloor+1}^{n} \bigg(h_{\mu^{\stackrel{0}{{\scriptscriptstyle 0}}}_{1,n+1}}(R_{j+1} \ | \ R_{j}) & \ - \ h_{\mu^{\stackrel{0}{{\scriptscriptstyle 0}}}_{1,n}}(R_j \ | \ R_{j-1}) \bigg)\\
\label{line3} + \ \sum_{k=1}^{\lfloor \frac{n}{2} \rfloor-1} \bigg(h_{\mu^{\stackrel{0}{{\scriptscriptstyle 0}}}_{1,n+1}}(R_{k} \ | \ R_{k+1}) & \ - \ h_{\mu^{\stackrel{0}{{\scriptscriptstyle 0}}}_{1,n}}(R_k \ | \ R_{k+1}) \bigg)\\
\label{line4} + \ h_{\mu^{\stackrel{0}{{\scriptscriptstyle 0}}}_{1,n+1}}\big(R_{\lfloor \frac{n}{2} \rfloor + 1} & \ | \ R_{\lfloor \frac{n}{2} \rfloor}\big).
\end{align}
Theorem~\ref{C1.5} implies that the measures $\mu^{\stackrel{0}{{\scriptscriptstyle 0}}}_{1,n+1}|_{R_{\lfloor \frac{n}{2} \rfloor}}$ and $\mu^{\stackrel{0}{{\scriptscriptstyle 0}}}_{1,n}|_{R_{\lfloor \frac{n}{2} \rfloor}}$ are ergodic. Then by Theorem~\ref{closerows} and Theorem~\ref{dbarholder}, it is clear that (\ref{line1}) is exponentially small in $n$, i.e. there exist constants $Q$ and $R$ independent of $n$ so that (\ref{line1}) $< Qe^{-Rn}$. We may rewrite any term in the sum (\ref{line2}) by Proposition~\ref{P2}:
\[
h_{\mu^{\stackrel{0}{{\scriptscriptstyle 0}}}_{1,n+1}}(R_{j+1} \ | \ R_{j}) - h_{\mu^{\stackrel{0}{{\scriptscriptstyle 0}}}_{1,n}}(R_j \ | \ R_{j-1})
\]
\[
= \left(h_{\mu^{\stackrel{0}{{\scriptscriptstyle 0}}}_{1,n+1}}(R_{j+1} \cup R_{j}) - h_{\mu^{\stackrel{0}{{\scriptscriptstyle 0}}}_{1,n}}(R_j \cup R_{j-1})\right) - \left(h_{\mu^{\stackrel{0}{{\scriptscriptstyle 0}}}_{1,n+1}}(R_{j}) - h_{\mu^{\stackrel{0}{{\scriptscriptstyle 0}}}_{1,n}}(R_{j-1})\right).
\]
By Theorem~\ref{closerows}, 
\[
\dbar\Big(\mu^{\stackrel{0}{{\scriptscriptstyle 0}}}_{1,n+1}|_{R_{j}}, \mu^{\stackrel{0}{{\scriptscriptstyle 0}}}_{1,n}|_{R_{j-1}}\Big) \leq Ae^{-Bj} \textrm{ and } \dbar\Big(\mu^{\stackrel{0}{{\scriptscriptstyle 0}}}_{1,n+1}|_{R_{j+1} \cup R_{j}}, \mu^{\stackrel{0}{{\scriptscriptstyle 0}}}_{1,n}|_{R_j \cup R_{j-1}}\Big) \leq Ae^{-Bj}.
\]
Since $j > \lfloor \frac{n}{2} \rfloor$ and all of the relevant measures are ergodic, (\ref{line2}) is exponentially small in $n$ by Theorem~\ref{dbarholder}. The proof that (\ref{line3}) is also exponentially small in $n$ is similar.

All that remains is to show that the leftover term (\ref{line4}) approaches $h$ at rate which is at least exponential in $n$. It suffices to show that (\ref{line4}) approaches any limit at all with rate at least exponential in $n$; by Lemma~\ref{L1}, $h_{n+1} - h_n$ approaches $h$ in the Ces\`{a}ro limit, and $h_{n+1} - h_n$ differs from (\ref{line4}) by an exponentially small amount. So, if (\ref{line4}) approaches a limit at all, it must be $h$.

We note that for any $n$, $\lfloor \frac{n+1}{2} \rfloor$ is either equal to $\lfloor \frac{n}{2} \rfloor$ or $\lfloor \frac{n}{2} \rfloor + 1$.
But by Theorem~\ref{closerows}, in either event, 
\[
\dbar\Big(\mu^{\stackrel{0}{{\scriptscriptstyle 0}}}_{1,n+2}|_{R_{\lfloor \frac{n+1}{2} \rfloor} \cup R_{\lfloor \frac{n+1}{2} \rfloor + 1}}, \mu^{\stackrel{0}{{\scriptscriptstyle 0}}}_{1,n+1}|_{R_{\lfloor \frac{n}{2} \rfloor} \cup R_{\lfloor \frac{n}{2} \rfloor} + 1}\Big) < Ae^{-B(\frac{n}{2})}.
\] 
\noindent
Then $h_{\mu^{\stackrel{0}{{\scriptscriptstyle 0}}}_{1,n+2}}\big(R_{\lfloor \frac{n+1}{2} \rfloor + 1} \ | \ R_{\lfloor \frac{n+1}{2} \rfloor}\big) - h_{\mu^{\stackrel{0}{{\scriptscriptstyle 0}}}_{1,n+1}}\big(R_{\lfloor \frac{n}{2} \rfloor + 1} \ | \ R_{\lfloor \frac{n}{2} \rfloor}\big)$ is exponentially small in $n$ by ergodicity and Theorem~\ref{dbarholder}, implying that (\ref{line4}) is exponentially Cauchy, therefore it approaches a limit with rate at least exponential in $n$, and we are done.

\begin{flushright}
$\blacksquare$\\
\end{flushright}

One application of Theorem~\ref{mainresult} is to the computability of the real number $h$. We first need to define our notion of computability.

\begin{definition}
A real number $\alpha$ is \textbf{computable in time} $f(n)$ if there exists a Turing machine which, on input $n$, outputs a pair $(p_n, q_n)$ of integers such that $|\frac{p_n}{q_n} - \alpha| < \frac{1}{n}$, and if this procedure takes less than $f(n)$ operations for every $n$. We say that $\alpha$ is \textbf{computable} if it is computable in time $f(n)$ for some function $f(n)$.
\end{definition}

Informally speaking, a real number $\alpha$ is computable if it is possible to give a finite description of $\alpha$ which allows someone to reconstruct as many digits of the decimal expansion of $\alpha$ as desired. For instance, $e$ is computable since we can describe it as the sum of the reciprocals of the factorials of nonnegative numbers. All algebraic numbers are computable, but there are many more computable numbers than algebraic (though still only countably many.) For an introduction to computability theory, see \cite{Ko}. 

\begin{theorem}\label{polycomp}
$h$ is computable in polynomial time. (There exists a polynomial $p(n)$ for which $h$ is computable in time $p(n)$.)
\end{theorem}

\noindent
\textit{Proof.} Recall from Section~\ref{intro} that for any $\mathbb{Z}$ nearest neighbor SFT $X$, $h^{top}(X)$ is the logarithm of the Perron eigenvalue of an associated matrix called its transition matrix. Since we will need a few relevant properties of these matrices, we quickly define them for $\mathbb{Z}$ nearest neighbor SFTs. Given a $\mathbb{Z}$ nearest neighbor SFT, which we assume without loss of generality to have alphabet $[1,|A|]$,
the transition matrix $B$ is a square $0$-$1$ matrix with size $|A|$, where $b_{ij}$ is $0$ if the adjacency $ij$ is not allowed and $1$ if the adjacency $ij$ is allowed.

Define, for any $n$, $B_n$ to be the transition matrix for $H_n = H_n(\mathcal{H})$. Then $B_n$ is a square matrix with size $s_n := LA_{\{1\} \times [1,n]}$. Since the horizontal adjacency conditions for $\mathcal{H}$ are symmetric ($ij$ is legal if and only if $ji$ is legal), the same is true for $H_n$, and so all $B_n$ are symmetric. For $\mathcal{H}$, the algorithm from \cite{Pi} mentioned in Section~\ref{intro} for generating any $B_n$ takes exponential time in $n$. (Briefly, one constructs $B_{n+1}$ from $B_n$ by arranging four copies of $B_n$ in a square, and then by replacing the right half of the upper-right copy of $B_n$, the upper half of the lower-left copy of $B_n$, and the entire lower-right copy of $B_n$ by $0$s. The number of operations taken to generate this matrix is of the same order as the number of operations it takes to write down the entries, of which there are exponentially many in $n$.) Also, $B_n$ is nonnegative real and symmetric, therefore it has all real eigenvalues, which we denote by $\lambda_{n,1}$, $\lambda_{n,2}$, \ldots, $\lambda_{n,s_n}$, where $\lambda_{n,1} \geq |\lambda_{n,2}| \geq \ldots \geq |\lambda_{n,s_n}|$. For any positive integer $k$, $\textrm{tr}((B_n)^k) = \sum_{i=1}^{s_n} \lambda_{n,i}^k$, and so if we assume $k$ to be even, then
\[
\lambda_{n,1}^k \leq \textrm{tr}((B_n)^k) \leq s_n \lambda_{n,1}^k.
\]
Since $s_n \leq 2^n$,
\[
\lambda_{n,1} \leq [\textrm{tr}((B_n)^k)]^{\frac{1}{k}} \leq 2^{\frac{n}{k}} \lambda_{n,1}.
\]
If we choose $k = 8^n$, then $k \geq n 4^n$, and so
\[
\lambda_{n,1} \leq [\textrm{tr}((B_n)^{8^n})]^{\frac{1}{8^n}} \leq 2^{4^{-n}} \lambda_{n,1}.
\]
Since $B_n$ is a $0$-$1$ matrix, $\lambda_{n,1}$ is less than or equal to the size $s_n$ of $B_n$, which is in turn less than $2^n$. Combining this with the fact that $2^{4^{-n}} \leq 1 + 4^{-n}$ yields $|\lambda_{n,1} - [\textrm{tr}((B_n)^{8^n})]^{\frac{1}{8^n}}| < 2^{-n}$. Also, the calculation of $[\textrm{tr}((B_n)^{8^n})]^{\frac{1}{8^n}}$ takes exponentially many steps in $n$; one simply needs to start with $B_n$ and square $3n$ times, then add the diagonal entries and take the result to the $\frac{1}{8^n}$ power. \\

Therefore, by investing exponentially many steps in $n$, it is possible to achieve approximations $\widetilde{h_{n+1}}$ and $\widetilde{h_n}$ which are exponentially close to $h_{n+1}$ and $h_n$ respectively, and then by Theorem~\ref{mainresult}, $\widetilde{h_{n+1}} - \widetilde{h_n}$ is exponentially close to $h$. 

In other words, there exist $C$, $D$, $E$, and $F$ so that for every $n$, there is an approximation, computable in less than $C e^{Dn}$ steps, which is within $E e^{-Fn}$ of $h$. But then for any integer $m$, $E e^{-F(n+1)} \leq \frac{1}{m} \leq E e^{-Fn}$ for some $n$, and so one can approximate $h$ to within $\frac{1}{m}$ in at most $C e^{D(n+1)}$ steps. Since $m \geq \frac{1}{E} e^{Fn}$, the number of steps required for the approximation is at most $Ce^D (mE)^{\frac{D}{F}}$, which is clearly a polynomial in $m$.

\begin{flushright}
$\blacksquare$\\
\end{flushright}

The fact that $h$ is computable follows from a more general result in \cite{HoMe}, but to our knowledge, very little was known about the rate. Another consequence of \cite{HoMe} is that there exist $\mathbb{Z}^2$ SFTs whose entropies are computable with arbitrarily poor time (along with entropies which are not computable at all!), so Theorem~\ref{polycomp} at least implies that $h$ is ``nice'' within the class of entropies of SFTs. Though not as good as a closed form, this is still satisfying; since $\mathcal{H}$ is the simplest possible nondegenerate $\mathbb{Z}^2$ SFT, one would hope for its entropy to be a relatively simple number.

\section{A counterexample}
\label{counterexample}

Interestingly, it is not true for all $\mathbb{Z}^2$ SFTs that $h_{n+1}(X) - h_n(X)$ converges to a limit. This was shown by an example in \cite{Pi}. However, this example was somewhat degenerate in that it was periodic, and in particular not topologically mixing.

\begin{definition}
A $\mathbb{Z}^d$ subshift $X$ is \textbf{topologically mixing} if for any finite rectangular prisms $S,T \subset \mathbb{Z}^d$, there exists $R_{S,T}$ so that for any translations $S'$ and $T'$ of $S$ and $T$ respectively such that $\|s'-t'\|_{\infty} > R_{S,T}$ for all $s' \in S'$ and $t' \in T'$, and for any globally admissible configurations $u \in L_{S'}(X)$ and $v \in L_{T'}(X)$, there exists $x \in X$ such that $x|_{S'} = u$ and $x|_{T'} = v$.
\end{definition}

In other words, $X$ is topologically mixing if it is possible to see any two globally admissible configurations at any desired locations within the same point of $X$, provided that you allow enough distance between them. Topological mixing is a strong condition for $\mathbb{Z}$ SFTs, and is a sufficient hypothesis for many theorems. However, for $\mathbb{Z}^d$ SFTs with $d>1$, topological mixing is a somewhat weak property. For many theorems in $\mathbb{Z}^d$ symbolic dynamics (see \cite{desai1}, \cite{desai2}, \cite{JM}, \cite{lightwood1}, and \cite{lightwood2}), it is necessary to assume a uniform mixing property, i.e. one where the distance required to see two globally admissible configurations simultaneously is independent of their size. 
There is a hierarchy of uniform mixing conditions in $\mathbb{Z}^d$, including block gluing, corner gluing, the uniform filling property, strong irreducibility, and square filling mixing. (See \cite{BPS} for definitions of and some exposition on the conditions in this hierarchy.)

We can modify the example from \cite{Pi} to see that the weakest uniform mixing condition, block gluing, is not enough to ensure convergence of $h_{n+1}(X) - h_n(X)$. 

\begin{definition}
A $\mathbb{Z}^d$ subshift $X$ is \textbf{block gluing} if there exists $R$ such that for any finite rectangular prisms $S,T \subset \mathbb{Z}^d$ satisfying  $\|s - t\|_{\infty} > R$ for all $s \in S$ and all $t \in T$, and for any globally admissible configurations $u \in L_S(X)$ and $v \in L_T(X)$, there exists $x \in X$ such that $x|_S = u$ and $x|_{T} = v$. We call the minimum such $R$ the \textbf{filling length} of $X$.
\end{definition}

We will not define any other uniform mixing conditions except to say that the stronger conditions have the same spirit, but enlarge the class of configurations which are considered. For instance, strong irreducibility is defined by considering any pair of globally admissible configurations, whether their shapes are rectangular prisms or something more complicated.

\begin{theorem}\label{Thmcounterexample}
There exists a block gluing $\mathbb{Z}^2$ nearest neighbor SFT $Y$ for which $\lim_{n \rightarrow \infty} h_{n+1}(Y) - h_n(Y)$ does not exist.
\end{theorem}

\begin{figure}[h]
\centering
\includegraphics[scale=0.4]{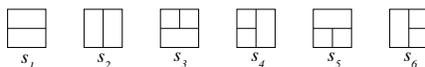}
\caption{Grid symbols in the alphabet $A$}
\label{F2}
\end{figure}

We begin by defining $Y$, which is a slightly different version of the SFT $X_{MS}^{(N)}$ defined in \cite{BPS}. The alphabet $A$ of $Y$ consists of the integers $0,1,\ldots,k$ for any $k > (8 \cdot 48^2)^2 = 339738624$, along with the symbols $s_1,s_2,s_3,s_4,s_5,s_6$ (illustrated in Figure~\ref{F2}), which we call grid symbols. The legal adjacent pairs of grid symbols are those where the line segments which meet the edges ``match up,'' and which do not yield parallel line segments at a unit distance which do not meet. For instance, $s_1 s_2$ is forbidden since the horizontal line segment meeting the right edge of $s_1$ does not match up with any horizontal line segment meeting the left edge of $s_2$, and the pairs $s_2 s_6$ and $s_3 s_4$ are forbidden since each pair would contain a pair of vertical line segments at unit distance which do not meet. Adjacencies between integers are as follows: $0$ may only appear horizontally adjacent to $0$, $0$ may not appear vertically adjacent to $0$, and non-$0$ integers may not be vertically adjacent. The only integer allowed to appear above a grid symbol is $0$, and there are no other restrictions on adjacencies between grid symbols and non-grid symbols.

The net effect of all of this is that any point $y \in Y$ has grid symbols partitioning the plane into rectangles (possibly infinite), and on each rectangle $y$ is labeled with integers, where the rows alternate between rows of all $0$s and rows of arbitrary strings of non-$0$ integers between $1$ and $k$. In any such rectangle (finite or infinite) with a bottom row, this row must be labeled with $0$s.

\

\begin{figure}[h]
\centering
\includegraphics[scale=0.4]{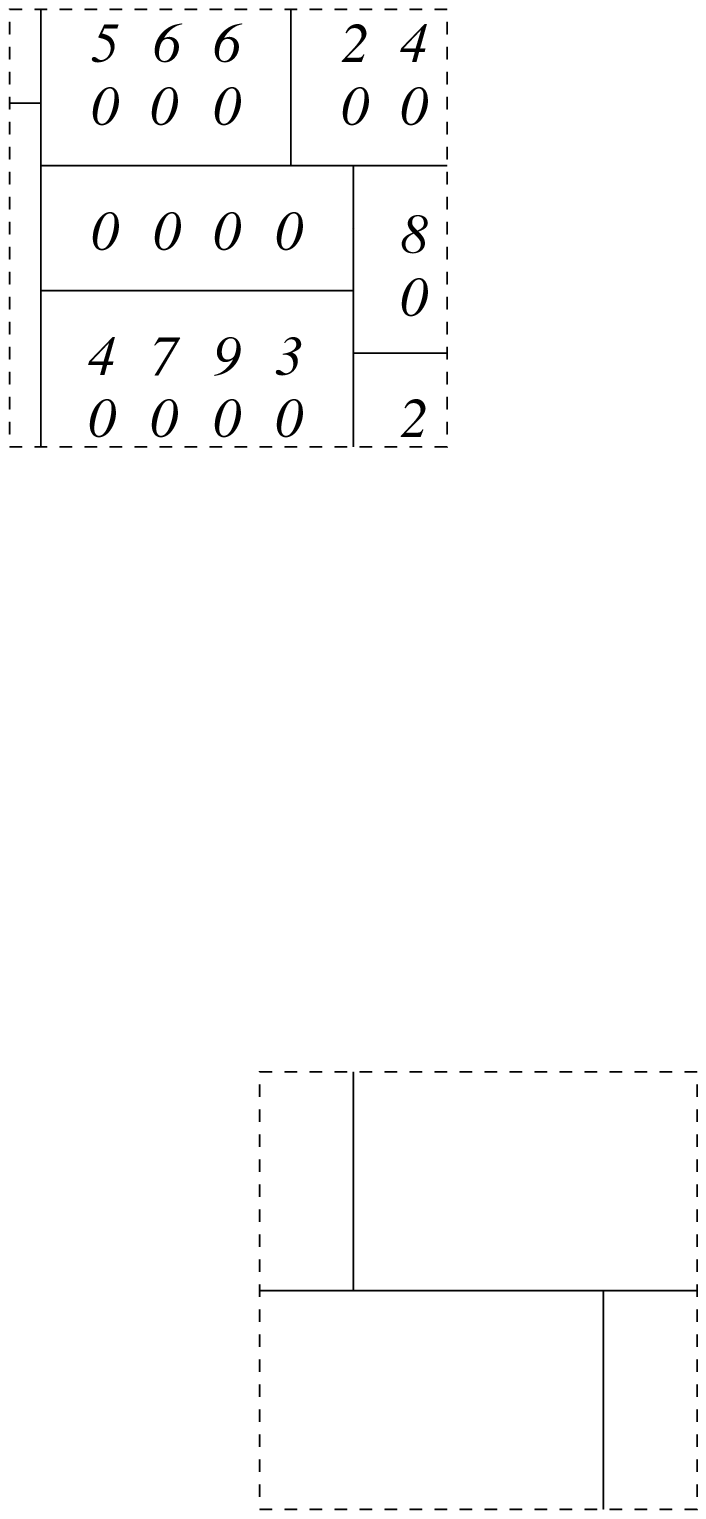}
\caption{$w$ and $w'$}
\label{F3}
\end{figure}

First, we will verify that $Y$ is block gluing with filling length $9$. Consider any two rectangular configurations $w$ and $w'$ which are globally admissible in $Y$. Without loss of generality, we assume that both $w$ and $w'$ have shape $[1,n]^2$. For any $v \in \mathbb{Z}^2$ with $\|v\|_{\infty} > n+9$, we will construct $x \in A^{\mathbb{Z}^2}$ for which $x|_{[1,n]^2} = w$ and $x|_{[1,n]^2 + v} = w'$. First, place $w$ and $w'$ at the corresponding locations, as in Figure~\ref{F3}. Clearly either the horizontal separation or vertical separation between $w$ and $w'$ is greater than $9$, and for now we assume that it is the vertical separation.

\begin{figure}[h]
\centering
\includegraphics[scale=0.4]{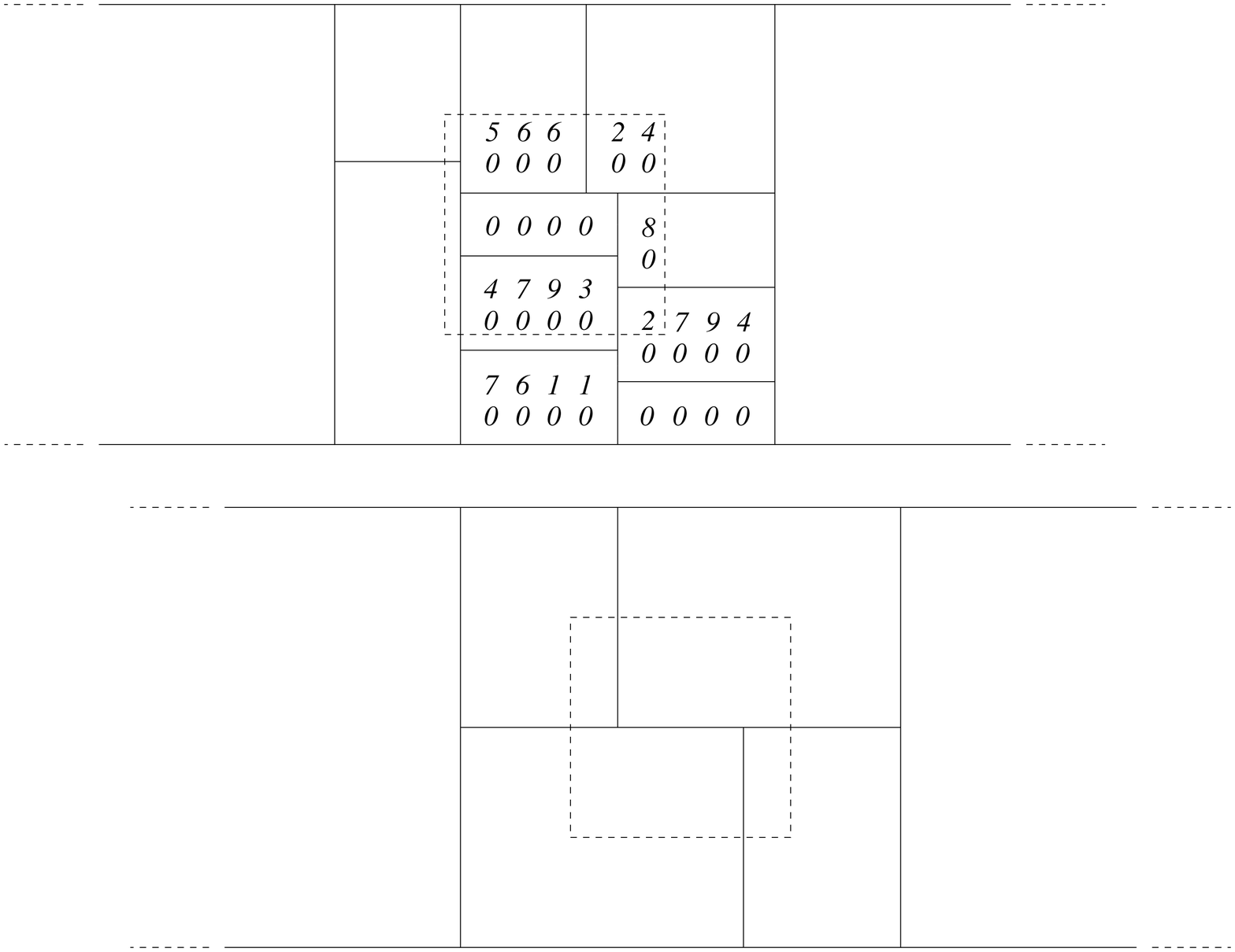}
\caption{Interpolating between $w$ and $w'$}
\label{F4}
\end{figure}

First, we will extend each of $w$ and $w'$ to a slightly larger square. We describe the procedure only for $w$, as the corresponding procedure for $w'$ is completely analogous. We begin by placing grid symbols on the border of $[-3,n+4]^2$, i.e. at a distance of $4$ from $w$. The top and bottom edges are labeled with horizontal lines (the symbol $s_1$), the left and right edges are labeled with vertical lines (the symbol $s_2$), the lower two corners are labeled with the symbol $s_3$, and the upper two corners are labeled with the symbol $s_5$. Denote the square $[-3,n+4]^2$ by $B$. For each edge of $w$, look for any grid symbols which contain a line segment which hits the boundary of $w$, and extend such segments to the corresponding edge of $B$ by using a string of grid symbols $s_1$ or $s_2$ (along with the proper ``joining'' symbol $s_3$, $s_4$, $s_5$, or $s_6$ when this string hits the edge.) This partitions $B$ into rectangles, which we would like to fill with integers. Any empty rectangles are easy to fill, and we can almost just complete the rectangles which already contain some integers from $w$ in a locally admissible way, but there is one slight problem; when continuing the pattern of alternating rows of $0$s and rows of non-$0$ integers begun by a partially filled rectangle, we could end up with a non-$0$ integer above one of the horizontal line grid symbols along the bottom edge of $B$, which is illegal. This is easily addressed though: before filling in any rectangles, consider any interval of integers on the bottom edge of $w$. If such an interval is made up of $0$s, place a horizontal line of grid symbols below it to end its rectangle. If an interval is made up of non-$0$ integers, place a row of $0$s below it, and then place a horizontal line of grid symbols below that. Again extend any incomplete paths to the boundary of $B$, and since each rectangle intersecting the bottom edge of $B$ is now empty, it is possible to fill all rectangles with integers, without changing $w$, in a locally admissible way. The resulting configuration on $B$ (and the corresponding one on $B'$) is locally admissible. To fill the rest of $\mathbb{Z}^2$, we simply extend the segments of horizontal line grid symbols on the top and bottom edges of both $B$ and $B'$ infinitely to the left and right, and fill in the resulting empty infinite rectangles with integers in any locally admissible way. (This procedure is illustrated in Figure~\ref{F4}.)  

If instead the horizontal separation between $w$ and $w'$ was at least $9$, then the only changes to the above construction would be to use grid symbols $s_4$ and $s_6$ on the corners of $B$ and $B'$, and to extend the segments of vertical line grid symbols on the left and right edges of $B$ and $B'$ infinitely upwards and downwards instead. Since $w$ and $w'$ must have been in one of these two situations, we have proved that $Y$ is block gluing.

\

We will now verify that $\lim_{n \rightarrow \infty} h_{n+1}(Y) - h_n(Y)$ does not exist. The basic idea is that most of the entropy is contributed only by the integer symbols in $A$, and that the entropy contributed by these symbols grows a lot when transitioning from a strip of height $2n$ to a strip of height $2n+1$, and not as much when transitioning from a strip of height $2n-1$ to a strip of height $2n$. 

Fix any $n,m \in \mathbb{N}$. We will bound $|LA_{[1,m] \times [1,2n-1]}(Y)|$ from above and below. The lower bound is easy: by considering configurations labeled by alternating rows of $0$s and non-$0$ integers with non-$0$ rows on the top and bottom, we quickly see that $|LA_{[1,m] \times [1,2n-1]}(Y)| \geq k^{mn}$. For the upper bound, we have to work a bit harder. Consider any configuration $w \in LA_{[1,m] \times [1,2n-1]}(Y)$ which contains $g$ grid symbols for some $g \in [1,(2n-1)m]$. We first bound from above the number of ways that these grid symbols can be placed. 

The key point in our argument is that since, in points of $Y$, each grid symbol must be adjacent to at least two other grid symbols and there are no grid symbols consisting of only a corner, there are no locally admissible closed finite loops of grid symbols. Therefore, any grid symbol in the interior of a locally admissible configuration with shape $[1,m] \times [1,2n-1]$ is part of a path of adjacent grid symbols which hits the boundary of $[1,m] \times [1,2n-1]$ at least twice (once entering, once leaving.) This enables us to design an algorithm which allows a Turing machine to recreate any locally admissible configuration of $g$ grid symbols on $[1,m] \times [1,2n-1]$ given a specific piece of input consisting of a finite ordered list $L$ of coordinates on the boundary of $[1,m] \times [1,2n-1]$ and a $g$-tuple ${\mathcal I}$ of instructions taken from a set of $48$ different commands. The list $L$ consists of sites on the border of the rectangle $[1,m] \times [1,2n-1]$. Every instruction in $\mathcal I$ is itself a $3$-tuple $(a_i,f_i,d_i)$ ($1\leq i\leq m$), where $a_i\in\{s_1,s_2,s_3,s_4,s_5,s_6\}$ represents one of the six grid symbols in the alphabet of $Y$, $f_i\in\{0,1\}$ is a flag that signals either ``revert'' or ``continue,'' and $d_i\in\{\text{up},\text{down},\text{right},\text{left}\}$ is one of the four standard directions in $\mathbb{Z}^2$. Now the Turing machine processes its input and puts down grid symbols on $[1,m] \times [1,2n-1]$ as follows: the machine starts by moving its writing-head to the coordinate given by the first entry in the list $L$ (if $L$ is empty, the algorithm stops here.) There it puts down the symbol $a_1$ from the first instruction in $\mathcal I$ starting a finite part of some path. If $f_1$ is ``continue,'' it moves its writing-head one step in the direction given by $d_1$, where it executes the next instruction in the same manner. If some $f_i$ is ``revert,'' the machine moves back along the grid symbols written so far until it comes to the first junction (one of the symbols $\{s_3,s_4,s_5,s_6\}$) where one of the three branches is a dead-end (i.e. the branch points to a place still inside $[1,m] \times [1,2n-1]$ where the machine has not already placed another grid symbol.) From there, the machine moves one step in the direction specified by $d_i$ and continues with the $(i+1)$th instruction. If there is no dead-end, the machine moves its writing-head to the next coordinate from the list $L$, where it starts another path of grid symbols using the next instruction from $\mathcal I$. After executing all commands in $\mathcal I$, the machine has placed exactly $g$ non-blanks. 

We claim that every locally admissible configuration $w$ consisting of $g$ grid symbols can be created by our Turing machine using some input. If $g=0$, clearly the empty input suffices. If $g>0$, then there is some grid symbol on the border of $[1,m] \times [1,2n-1]$, which we can take to be the first site in $L$. Then, follow any path of adjacent grid symbols in $w$, recording the proper entries of $\mathcal{I}$, until you either run into the border of $[1,m] \times [1,2n-1]$, or will be forced to run into an already visited grid symbol. If you have visited all $g$ grid symbols in $w$, then you are done. Since $w$ does not contain closed finite loops, if there are still unvisited grid symbols in $w$, then they are all either connected to an already visited grid symbol or connected to the border of $[1,m] \times [1,2n-1]$ by a path of adjacent grid symbols. So, we can record an entry of $\mathcal{I}$ with $f_i$ ``revert,'' and either move back to the first place along your path where you could continue to unvisited grid symbols, or, if this is impossible, begin with an unvisited grid symbol on the border of $[1,m] \times [1,2n-1]$, append this site to $L$, and continue. In this fashion, we can eventually visit all $g$ grid symbols in $w$, simultaneously recording the input which will recreate $w$.

Therefore, the number of different input ``programs'' gives an upper bound on the number of ways to place $g$ grid symbols on $[1,m] \times [1,2n-1]$ in a locally admissible way. By overestimating the number of lists $L$ by the number of subsets of the boundary of $[1,m] \times [1,2n-1]$, we get an upper bound of $2^{4n + 2m - 6} 48^g$.

Now, fix any locally admissible assignment of $g$ grid symbols. We wish to bound from above the number of ways to fill in the leftover rectangles with integers in a locally admissible way. For any $w \in LA_{[1,m] \times [1,2n-1]}(Y)$, consider a column of $w$ which has $h$ grid symbols in it. This column consists of alternating intervals of integers and grid symbols. Due to the restriction that non-$0$ integers cannot appear above grid symbols, each one of these intervals of integers has at most half non-$0$ integers, except possibly for the bottom-most interval, which could have one more non-$0$ integer than $0$. This means that the total number of non-$0$ integers in the column is at most $n - \frac{h}{2}$. Since the only choice for each interval of integers is whether its bottom-most integer is $0$ or non-$0$ and which non-$0$ integers to use, and since only the bottom-most interval admits a choice about whether its bottom-most integer is $0$ or non-$0$, this implies that the total number of ways of filling the leftover portion of this column with integers is at most $2 \cdot k^{n-\frac{h}{2}}$. Therefore, the total number of ways to extend any fixed locally admissible grid symbol configuration containing $g$ grid symbols to a locally admissible configuration on all of $[1,m] \times [1,2n-1]$ is at most $2^m k^{mn - \frac{g}{2}}$, and so
\[
|LA_{[1,m] \times [1,2n-1]}(Y)| \leq \sum_{g=0}^{(2n-1)m} 2^{4n+2m-6} 48^g 2^m k^{mn - \frac{g}{2}} \leq ((2n-1)m+1)2^{4n + 3m - 6} k^{mn}.
\]
(Here the last inequality uses the fact that $k > 48^2$.) Combining with the earlier lower bound on $|LA_{[1,m] \times [1,2n-1]}(Y)|$, taking logarithms, dividing by $m$, and letting $m \rightarrow \infty$ yields the bounds $n \ln k \leq h_{2n-1}(Y) \leq n \ln k + \ln 8$.

We will now achieve similar bounds on $|LA_{[1,m] \times [1,2n]}(Y)|$. Again, we may arrive at a lower bound by considering only configurations of integers: $|LA_{[1,m] \times [1,2n]}(Y)| \geq k^{mn}$. By the same proof as before, the number of ways that $g$ grid symbols can be placed on $[1,m] \times [1,2n]$ in a locally admissible way is less than $2^{4n + 2m - 4} 48^g$. Also by the same proof, the number of ways to fill a column with $h$ grid symbols in a locally admissible way is at most $2 \cdot k^{n - \frac{h-1}{2}}$. We note that if $h=0$, then the number of ways to fill the column is clearly $2k^n$, and so our upper bound is $\min(2k^n, 2 \cdot k^{n - \frac{h-1}{2}})$. The number of ways to complete a fixed locally admissible grid symbol configuration containing $g$ grid symbols to a locally admissible configuration on all of $[1,m] \times [1,2n]$ is then at most $\min(2^m k^{mn}, 2^m k^{mn - \frac{g-m}{2}})$, and we get the upper bound
\begin{multline}\notag
|LA_{[1,m] \times [1,2n]}(Y)| \leq \sum_{g = 0}^{2nm} 2^{4n + 2m - 4} 48^g 2^{m} \min(k^{mn}, k^{mn - \frac{g-m}{2}})\\
= 2^{4n + 3m - 4} \sum_{g = 0}^{2mn} 48^g \min(k^{mn}, k^{mn - \frac{g-m}{2}})
\end{multline}
\[
\leq 2^{4n + 3m - 4} \Big[ \sum_{g = 0}^{2m} 48^{2m} k^{mn} + \sum_{g = 2m + 1}^{2mn} 48^g k^{mn - \frac{g}{4}}\Big] \leq 2^{4n + 3m - 4} (2mn + 1) 48^{2m} k^{mn}.
\]

(The last inequality uses the fact that $k > 48^4$.) Combining with the earlier lower bound on $|LA_{[1,m] \times [1,2n]}(Y)|$, taking logarithms, dividing by $m$, and letting $m \rightarrow \infty$ yields the bounds $n \ln k \leq h_{2n}(Y) \leq n \ln k + \ln(8 \cdot 48^2)$. But then for any $n$, $h_{2n+1}(Y) - h_{2n}(Y) \geq \ln k - \ln (8 \cdot 48^2)$ and $h_{2n}(Y) - h_{2n-1}(Y) \leq \ln(8 \cdot 48^2)$. Since $k > (8 \cdot 48^2)^2$, this means that there exists $\epsilon > 0$ so that $h_{2n+1}(Y) - h_{2n}(Y) > h_{2n}(Y) - h_{2n-1}(Y) + \epsilon$ for all $n$, and so $h_{n+1}(Y) - h_n(Y)$ does not approach a limit as $n \rightarrow \infty$.

\begin{flushright}
$\blacksquare$\\
\end{flushright}

\section{Questions}
\label{questions}

There are several questions which suggest themselves from this work. Firstly, though we have shown that $h_{n+1} - h_n \rightarrow h$ at a rate which is at least exponential, we have not been able to give any explicit bound for this rate. 

\begin{question}\label{rate}
Is it possible to give explicit values of $A$ and $B$ for which $h_{n+1} - h_n < Ae^{-Bn}$?
\end{question}

The answer to this question would be interesting both because it might allow us to improve the known bounds on $h$ and also because it would allow us to give an explicit polynomial upper bound on the time of computability of $h$. In order to find such $A$ and $B$, it would be sufficient to give an explicit such $A$ and $B$ for $p=0.5$ in Theorem~\ref{expdecay}, but it seems that finding these is somewhat difficult. We note that for much smaller percolation probabilities than $0.5$, giving explicit values for $A$ and $B$ is easy. For instance, if $p < 0.25$, then since there are less than $4^t$ paths from $0$ to $\partial([-n,n]^2)$ of length $t$ for any $t$, $P_p(0 \leftrightarrow \partial([-n,n]^2)) < \sum_{t = n}^{\infty} (4p)^t = \frac{1}{1-4p}(4p)^n$.

\begin{question}\label{ext}
Is it possible to extend these methods to a larger class of $\mathbb{Z}^2$ SFTs?
\end{question}

The difficulty here is that our proof relies on two important properties of $\mathcal{H}$. First, there must be some (possibly site-dependent) ordering on the alphabet for which the fundamental Theorem~\ref{Th1} is true, and this does not seem to be true for all shifts of finite type. Secondly, in order to use the methods of \cite{vdBS} to prove exponential closeness of the relevant measures with respect to $\dbar$, the SFT must satisfy a quite restrictive property related to conditional probability of disagreement at a pair of sites given their neighbors. (For most $\mathbb{Z}^2$ SFTs, the $0.5$ in Theorem~\ref{Th9} becomes a number larger than $p_c$, which means that we cannot show exponential decay.) So far, we have not been able to find any nondegenerate $\mathbb{Z}^2$ SFTs besides the hard square shift which have both of these properties, but it is possible that with a slightly different method, one could consider a wider class of systems.

\begin{question}\label{higherd}
Is it possible to apply these methods to the $\mathbb{Z}^d$ hard square shift for $d > 2$?
\end{question}

The difficulty here is that already $p_c(\mathbb{Z}^d) < 0.5$ for $d=3$ (\cite{CR}), which causes a problem with using Theorem~\ref{Th9} to imply exponential decay of $\dbar$ distance.

\section*{acknowledgements}
The author would like to thank the anonymous referee for useful corrections and comments, Larry Pierce for bringing this interesting problem to his attention and for several useful discussions, Lior Silberman for pointing out the relationship between Theorem~\ref{mainresult} and computability of $h$, and Joel Friedman for a helpful comment regarding computation of the approximations $\widetilde{h_n}$ of $h_n$ in Theorem~\ref{polycomp}. The author would also like to thank Brian Marcus for many extremely useful discussions, advice in the organization of this paper, and for pointing out a drastic simplification in the proof of Theorem~\ref{mainresult}.

\end{document}